\documentclass[reqno,11pt]{amsart}
\usepackage[margin=1in]{geometry}
\usepackage{amsfonts,amssymb,mathrsfs}
\usepackage[hidelinks]{hyperref}
\hypersetup{pdftitle={Nonlinear stability of viscous shock waves for the one--dimensional Navier--Stokes equations with singular viscosity},
	pdfauthor={Guanghan Zhao, Jingyu Li, Xiaowen Li, Changchun Liu, Ming Mei}}
\allowdisplaybreaks[4]
\numberwithin{equation}{section}
\newtheorem{theorem}{Theorem}[section]
\newtheorem{lemma}{Lemma}[section]
\newtheorem{proposition}{Proposition}[section]
\theoremstyle{remark}

\newcommand{\R}{\mathbb R}

\begin{document}
	\title[Stability of shocks with singular viscosity]
	{Nonlinear stability of viscous shock waves for the one--dimensional
		Navier--Stokes equations \\with singular viscosity}
	\author[G. Zhao, J. Li, X. Li, C. Liu, and M. Mei]
	{Guanghan Zhao, Jingyu Li, Xiaowen Li, Changchun Liu$^\dagger$, and Ming Mei$^*$}
	\address{Guanghan Zhao\hfill\break School of Mathematics, Jilin
		University, Changchun 130012, China}
	\email{zhaogh24@mails.jlu.edu.cn}
	\address{Jingyu Li\hfill\break School of Mathematics and Statistics,
		Northeast Normal University, Changchun 130024, China}
	\email{lijy645@nenu.edu.cn}
	\address{Xiaowen Li\hfill\break School of Mathematical Science,
		Peking University, Beijing 100871, China}
	\email{xwli@math.pku.edu.cn}
	\address{Changchun Liu\hfill\break School of Mathematics, Jilin
		University, Changchun 130012, China}
	\email{liucc@jlu.edu.cn}
	\address{Ming Mei\hfill\break School of Mathematics and Statistics,
		Jiangxi Normal University, Nanchang, Jiangxi 330022, China\hfill\break
		and Department of Mathematics, Champlain College Saint-Lambert,
		Quebec J4P 3P2, Canada\hfill\break
		and Department of Mathematics and Statistics, McGill University,
		Montreal, Quebec H3A 2K6, Canada}
	\email{ming.mei@mcgill.ca}
	
	\thanks{$^\dagger$Corresponding author.}
	\thanks{$^*$Corresponding author.}
	\date{}
	
	\begin{abstract}
		We study the nonlinear stability of viscous shock waves for the one--dimensional
		isentropic Navier--Stokes system with fast diffusion and singular viscosity. The waves connect a positive left state to zero specific volume at the right far field, where the density becomes unbounded and the viscous coefficient is singular. Under a vanishing--stress condition, we first
		establish the existence and uniqueness, up to a shift, of a monotone
		profile with an algebraic right tail. We then prove global existence and convergence to the profile for sufficiently small perturbations in weighted Sobolev spaces. The main difficulties stem from the singular nonlinear viscous term. We note that the existing approach for treating scalar conservation law with singular viscosity can not be applied to Navier--Stokes equations. We propose some new technique: after introducing an effective velocity and using an anti--derivative method, we establish unusual interactive--energy--estimates and thereby the expected energy estimates for velocity by coupling the two perturbation equations; furthermore, we introduce an auxiliary functional to eliminate uncontrollable terms. The artful weighted--energy estimates play a crucial role for treating singularity.
	\end{abstract}
	\maketitle
	
	\tableofcontents
	
	\section{Introduction}
	We consider the one--dimensional isentropic Navier--Stokes system in
	Lagrangian coordinates,
	\begin{equation}
		\begin{cases}
			v_t-u_x=0,\\
			u_t+p(v)_x=\mu\left(\dfrac{u_x}{v}\right)_x,
		\end{cases}
		\qquad x\in\R,\quad t>0,
		\label{1-1}
	\end{equation}
	with initial data
	\begin{equation}
		(v,u)(x,0)=(v_0,u_0)(x)\longrightarrow(v_\pm,u_\pm),
		\quad\text{as }x\to\pm\infty.
		\label{1-2}
	\end{equation}
	Here $v>0$ is the specific volume, $u$ is the velocity, $\mu>0$ is constant and $p$ is the pressure. We assume
	\begin{equation}
		v_->v_+=0.
		\label{v-}
	\end{equation}
	The singular end state
	\begin{equation*}
		v_+=0
	\end{equation*}
	corresponds to an unbounded density $\rho=1/v$ at the right far field.
	In particular, it is an infinite--density state, not a vacuum state.
	Throughout the paper, we impose
	\begin{equation}
		p\in C^4([0,\infty)),\quad p'(v)<0,\quad p''(v)>0,
		\quad \hbox{for all } v\geq0.
		\label{1-4}
	\end{equation}
	A typical example is $p(v)=(1+v)^{-\gamma}$, $\gamma>0$.
	These assumptions describe a generalized pressure law with a finite value
	at $v=0$; they exclude the usual polytropic law $p(v)=av^{-\gamma}$.
	
	\medskip
	
	{\bf Background of study.} Early studies of shock layers and fluid models provided useful background for viscous shock stability theory. Gilbarg \cite{GD} investigated the existence and limiting behavior of one--dimensional shock layers, while Kawashima, Matsumura and Nishida \cite{KSMANT} studied the fluid--dynamical approximation of the Boltzmann equation at the Navier--Stokes level. For the inviscid case, Nishida and Smoller \cite{NTSAJ} established global--solution results of certain nonlinear hyperbolic conservation laws. These works provide relevant physical background of \eqref{1-1}.
	
	When the system \eqref{1-1} is without singular viscosity, Matsumura and Nishihara \cite{MANK}
	established a classical nonlinear stability result under sufficiently small initial perturbations and a zero--mass condition. Goodman \cite{GJ} developed the related stability theory for scalar viscous conservation laws. Liu \cite{LT} further considered phase--shifts and diffusion waves generated by perturbations. General developments in viscous shock stability and physical--viscosity systems can also be found in \cite{LTN,LTZY}.
	
	While studying shock waves of compressible fluids, stability methods were developed for other one--dimensional systems. Matsumura and Nishihara \cite{MANKA} and Mei \cite{MMS} studied viscous scalar conservation laws with nonconvex fluxes. For the system of viscoelasticity, we refer to \cite{KSMA,MAMMN,MMNK,NK}. These results have contributed to the development of one--dimensional wave stability theory.
	
	Subsequent research was extended to more general systems, leading to the development of various methods for stability analysis. Mascia and Zumbrun \cite{MCZK} studied small--amplitude Lax shocks in symmetric hyperbolic--parabolic systems, including the compressible Navier--Stokes equations. In \cite{BBHJL,HJLOZK}, spectral stability of isentropic Navier--Stokes shock waves under large-amplitude conditions was investigated. Meanwhile, nonisentropic Navier--Stokes systems with large perturbations have also been considered; see \cite{HBTSZL,TSZL}. For the isentropic case, Wang, Zhao and Zou \cite{WTZHZQ} treated a class of large perturbations that allow initial density profiles with large oscillations.
	
	Density--dependent viscosity is particularly important when vacuum may occur. Matsumura and Wang \cite{MAWY} proved asymptotic stability of viscous shock waves under conditions of the viscosity law and a zero--mass assumption. Jiu, Wang and Xin \cite{JQWYXZ} studied large--time behavior near rarefaction waves whose far--field states include vacuum. Cui \cite{CM} considered the issue of vanishing viscosity in the presence of shock waves. He and Huang \cite{HLHF} proved nonlinear stability of large--amplitude shock profiles in the isentropic Navier--Stokes model with density--dependent viscosity. Related results on shock stability appear in \cite{HLTSWT,VAYL}.
	
	Boundary effects have also been studied. Matsumura and Mei \cite{MAMM} investigated the asymptotic
	behavior of viscous shock waves in the presence of boundaries. Liu and Cui \cite{LJCM} studied density--dependent viscosity shocks in a half--space, while Kang, Oh and Wang \cite{KMJHO} subsequently analyzed outflow shocks in the barotropic Navier--Stokes equations on a half--line.
	
	As the analysis moved beyond a single traveling wave, researchers also studied composite wave patterns. Huang and Matsumura \cite{HFMA} investigated the asymptotic stability of a composite wave consisting of two viscous shock waves. Kang, Vasseur and Wang \cite{KMJ} proved the long--time stability of a composite wave consisting of a viscous shock and a rarefaction wave, while Han, Kang and Kim \cite{HSKM} studied two viscous shock waves with independently small strengths.
	
	In addition to the energy method for proving stability, relative entropy method has also been applied to scalar or systems of conservation laws. Leger \cite{LN} studied shocks in scalar conservation laws, while Leger and Vasseur \cite{LNVA} treated shocks and contact discontinuities in systems of conservation laws under non--$BV$ perturbations. For scalar viscous conservation laws, Kang and Vasseur
	\cite{KMJVAL} obtained an $L^{2}$--type contraction result. Under large perturbations, Kang and Vasseur
	\cite{KMJVA} obtained a stability result of small barotropic Navier--Stokes shock waves.
	
	For scalar equations with singular diffusion, the stability of shock and rarefaction waves has been investigated in \cite{Fan-Mei-Xiang,Li-Li-Mei-Nave,Li-Li-Mei,Li-Mei,Xu-Mei-Nave-Sheng}. These works consider equations of the form
	\[
	u_t+f(u)_x=\mu\left(\frac{u_x}{u^m}\right)_x,\qquad m>0,
	\]
	and use singular weights adapted to the wave profile.
	
	\medskip
	
	{\bf Main difficulty and strategy.} Our aim is to establish the existence and nonlinear stability of viscous shocks connecting $v_->0$ to $v_+=0$. The profile of zero state at far field makes the coefficient $1/V$ unbounded. In this paper, we first establish the existence and uniqueness, up to translation, of a monotone viscous shock profile $(V,U)(x-st)$ connecting $(v_-,u_-)$ and $(0,u_+)$. We then prove the nonlinear asymptotic stability of $(V,U)$ under sufficiently small perturbations in suitably weighted Sobolev spaces.
	
	In  the previous studies, a standard strategy is to introduce the effective velocity
	\begin{equation*}
		b=u-\mu v_x/v,
	\end{equation*}
	first establish nonlinear stability of the transformed shock profile $(V,B)$ by the singular weighted--energy method, and then transfer the resulting stability estimate back to the original profile $(V,U)$. However, due to the presence of singularity in the system \eqref{1-1}, although the energy estimate provides a weighted space--time $L^{2}$ bound of the effective--velocity perturbation $\psi_{z}$, this bound is insufficient to close the subsequent $L^{2}$ estimate of $\phi_{z}/V$, where $\phi_{z}$ denotes the specific--volume perturbation.
	
	To overcome this difficulty, some new technique are expected. We directly study nonlinear stability of $(V,U)$. By exploiting the coupling between the two perturbation equations, we derive unusual interactive--energy--estimates to obtain a weighted space--time $L^{2}$ bound of the velocity perturbation $q$. To cancel the terms in this estimate that cannot be controlled directly, we introduce an auxiliary functional and establish the suitable estimate for it, thereby closing the desired energy estimate. Then, using the relationship between $\psi_{z}$ and $q$, we rewrite the perturbation system such that it involves only the perturbations $\phi_{z}$ and $q$. On this basis, we establish the corresponding nonlinear stability estimates of $(V,U)$ by artful weighted--energy method and the weighted space--time $L^{2}$ bound of $q$. Finally, the cut--off technique closes the argument.
	
	\medskip
	
	Before stating our main results, let us introduce some notations as follows.

	\noindent\textbf{Notation.}
	The letter $C$ denotes a positive constant whose value may vary from line to line. The integrals $\int_{\mathbb{R}} f(x) \mathrm{d}x$ and $\int_0^t \int_{\mathbb{R}} f(x, \tau) \mathrm{d}x \mathrm{d} \tau$ will be expressed as $\int f(x)$ and $\int_0^t \int f(x, \tau)$, respectively.
	We write $\|f\|=\|f\|_{L^2(\R)}$ and
	$\|f\|_k=\|f\|_{H^k(\R)}$ $(k \geq 0)$. For a positive weight $w$, set
	\[
	\|f\|_w^2=\int w|f|^2,\qquad
	\|f\|_{k,w}^2=\sum_{j=0}^k\int w|\partial_z^j f|^2.
	\]
	The corresponding spaces are denoted by $L_w^2$ and $H_w^k$.
	For a vector, $\| (f,g,h) \|^{2}$ is equivalent to $\| f \|^{2} + \| g \|^{2} + \| h \|^{2}$.
	The notation $f\sim g$ means $f/g\to1$ at the specified limit.
	
	\medskip
	
	{\bf Viscous shock waves.} An admissible viscous shock profile is a nonconstant pair
	$(V,U)\in C^4(\R)$, with $V>0$, such that
	\begin{equation}
		(v,u)(x,t)=(V,U)(z),\qquad z=x-st,
		\label{1-111}
	\end{equation}
	solves \eqref{1-1}. Substituting it into \eqref{1-1} yields
	\begin{equation}
		\begin{cases}
			-sV_z=U_z,\\
			-sU_z+p(V)_z=\mu\left(\dfrac{U_z}{V}\right)_z.
		\end{cases}
		\label{1-5}
	\end{equation}
	We require
	\begin{equation}
		\begin{gathered}
			V(-\infty)=v_-,\qquad V(+\infty)=0,\qquad U(\pm\infty)=u_\pm,\\
			V_z(\pm\infty)=U_z(\pm\infty)=0,\qquad
			\lim_{z\to\pm\infty}\frac{U_z(z)}{V(z)}=0.
		\end{gathered}
		\label{1-6}
	\end{equation}
	The last condition is a vanishing viscous--stress condition. Since $U_z\to0$ does not imply $U_z/V\to0$ when $V\to0$, this condition must be specified separately.
	Under \eqref{1-6}, the Rankine--Hugoniot conditions for a profile connecting $(v_{-},u_{-})$ to $(0,u_{+})$ are
	\begin{equation}
		u_+-u_-=sv_-,\qquad s^2=\frac{p(0)-p(v_-)}{v_-}.
		\label{1-7}
	\end{equation}
	The connecting profile necessarily has $s<0$. Strict convexity of $p$ gives the
	Lax's entropy condition
	\begin{equation}
		-\sqrt{-p'(0)}<s<-\sqrt{-p'(v_-)}.
		\label{1-8}
	\end{equation}
	
	\begin{theorem}[Existence and uniqueness of the profile]\label{theorem1.1}
		Assume that \eqref{1-4} holds. An admissible profile satisfying \eqref{1-6} exists
		if and only if \eqref{1-7} and \eqref{1-8} hold with $s<0$. The profile is unique up to a shift and satisfies
		\begin{equation}
			V_z<0,\qquad U_z<0,\qquad\text{on }\R.
			\label{1-9}
		\end{equation}
		Moreover, for some $A_->0$ depending on the translation,
		\begin{equation}
			V(z)\sim\frac{1}{a_+z}\quad(z\to+\infty),\qquad
			v_--V(z)\sim A_-e^{\lambda_-z}\quad(z\to-\infty),
			\label{1-911}
		\end{equation}
		where
		\[
		a_+=\frac{p'(0)+s^2}{s\mu}>0,\qquad
		\lambda_-=-\frac{v_-}{s\mu}\bigl(p'(v_-)+s^2\bigr)>0.
		\]
	\end{theorem}
	In the original system \eqref{1-1}, the nonlinear singular viscous term $(\frac{u_{x}}{v})_{x}$ makes the stability analysis of shock waves $(V,U)$ quite difficult. To handle the singular coefficient $\frac{1}{v}$ in the viscous term, we introduce the effective velocity
	\begin{equation}
		b=u-\mu\frac{v_x}{v}.
		\label{1-10}
	\end{equation}
	Substituting \eqref{1-10} into \eqref{1-1}, we obtain
	\begin{equation}
		\begin{cases}
			v_t-b_x=\mu\left(\dfrac{v_x}{v}\right)_x,\\
			b_t+p(v)_x=0,
		\end{cases}
		\label{1-11}
	\end{equation}
	with the initial data
	\begin{equation}
		(v,b)(x,0) = (v_0,b_0)(x) \to (v_\pm,b_\pm), \quad \text{as } x \to \pm \infty,
		\label{1-12}
	\end{equation}
	where
	\begin{equation*}
		b_0=u_0-\mu\frac{v_{0x}}{v_0}.
	\end{equation*}
	For an admissible profile, the stress condition in \eqref{1-6} gives $b_\pm=u_\pm$.
	
	Effective velocity formulations for related compressible models have been applied in
	\cite{HLHF, KMJVA, VAYL}; for viscoelastic systems, see also \cite{LJGZ}.
	The relevant structures of weak solutions also appear in \cite{BDBD,BDBDO}.
	
	Next, based on Theorem \ref{theorem1.1}, we state the existence of viscous shock waves for the system \eqref{1-11} for the effective velocity.
	
	\begin{theorem}\label{theorem1.2}
		Let $(V,U)$ be the profile in Theorem \ref{theorem1.1}. Define
		\[
		B=U-\mu\frac{V_z}{V}.
		\]
		Then the system \eqref{1-11} has a unique (up to a shift) smooth monotone traveling wave solution $(V,B)(x-st)$, with the end states
		$(v_-,u_-)$ and $(0,u_+)$, satisfying
		\begin{equation}
			V_z<0,\qquad B_z<0 \quad\text{on }\;\R.
			\label{1-14}
		\end{equation}
	\end{theorem}
	
	\medskip
	
	{\bf Stability of viscous shock waves.} Now we are going to state the stability of viscous shock waves. Let us denote
	\begin{equation}
		V_0(x)=V(x+x_0), \quad U_0(x)=U(x+x_0), \quad B_0(x)=B(x+x_0).
		\label{1-151}
	\end{equation}
	Assume that $v_0-V$ and $b_0-B$ are integrable. Using the end states of $V$ and $B$, we obtain
	\begin{equation}
		\begin{aligned}
			\int \binom{v_{0}(x)-V_{0}(x)}{b_{0}(x)-B_{0}(x)} \mathrm{d} x = & \int \binom{v_0(x)-V(x)}{b_0(x)-B(x)} \mathrm{d} x + \int \binom{V(x)-V_{0}(x)}{B(x)-B_{0}(x)} \mathrm{d} x
			\\
			= & \int \binom{v_0(x)-V(x)}{b_0(x)-B(x)} \mathrm{d} x-x_0\binom{0-v_{-}}{b_{+}-b_{-}}.
		\end{aligned}
		\label{3-2}
	\end{equation}
	We assume that there exists a shift $x_{0} \in \mathbb{R}$ such that
	\begin{equation}
		\int\binom{v_0(x)-V(x)}{b_0(x)-B(x)} \mathrm{d}x = x_0\binom{-v_-}{b_+-b_-}.
		\label{3-1}
	\end{equation}
	For this choice of $x_{0}$, \eqref{3-2} gives the zero--mass conditions
	\begin{equation*}
		\int \binom{v_{0}(x)-V_{0}(x)}{b_{0}(x)-B_{0}(x)} \mathrm{d} x = \binom{0}{0}.
	\end{equation*}
	Let us define
	\begin{equation*}
		(\phi_0,\psi_0)(x)=\int_{-\infty}^x
		(v_0(y)-V_0(y),\ b_0(y)-B_0(y))\,\mathrm{d}y.
	\end{equation*}
	Then $(\phi_0,\psi_0)(\pm\infty)=0$. The nonlinear stability of the shock wave $(V,U)(x-st)$ of \eqref{1-1}--\eqref{1-2} is stated in the following theorem.
	
	\begin{theorem}[Stability of viscous shock waves]\label{theorem1.3}
		Let $(V,B)(x-st+x_{0})$ be a pair of shock profiles obtained in Theorem \ref{theorem1.2}. Assume that $v_{0} - V_{0} \in L^{1}(\mathbb R)$, $b_{0} - B_{0} \in L^{1}(\mathbb R)$ and \eqref{3-1} holds. $V_{0}$, $U_{0}$ and $B_{0}$ are defined in \eqref{1-151}. Set $w_{2,0}=V_0^{-2}$. There exists $\delta>0$
		such that, if
		\begin{align*}
			\|\phi_0\| + \|\psi_0\|_2 + \|v_0-V_0\|_{1,w_{2,0}} + \|u_0-U_0\|_{1,w_{2,0}}\leq\delta,
		\end{align*}
		then \eqref{1-1}--\eqref{1-2} has a unique global solution $(v,u)(x,t)$ satisfying
		\[
		v-V \in C([0,\infty);H^1_{w_{2}}),\quad u-U \in C([0,\infty);H^1_{w_{2}}),
		\]
		where $w_{2} = V^{-2}$. Moreover, the solution converges uniformly to the shifted shock profile
		\begin{equation}
			\sup_{x\in\R}|(v,u)(x,t)-(V,U)(x-st+x_0)|\to 0,
			\quad\text{as }\;t\to\infty.
			\label{1-15}
		\end{equation}
	\end{theorem}
	
	The paper is organized as follows. In Section~\ref{sec:profiles}, we establish the existence, uniqueness, monotonicity and asymptotic behavior of $(V,U)$, and then of $(V,B)$. Section~\ref{sec:stability} is devoted to the weighted energy estimates and the proof of nonlinear stability. Appendix~\ref{sec:local} proves local existence.

	\section{Existence of shock profiles}\label{sec:profiles}
	\begin{proof}[Proof of Theorem \ref{theorem1.1}]
		Suppose first that an admissible profile $(V,U)(x-st)$ exists. Integrating the first
		equation of \eqref{1-5} gives
		\begin{equation}
			u_+-u_-=sv_-.
			\label{2-111}
		\end{equation}
		Since $u_{+} \neq u_{-}$, \eqref{2-111} gives $s \neq 0$. By $U_z=-sV_z$, the second equation of \eqref{1-5} becomes
		\begin{equation}
			s^2V_z+p(V)_z=-s\mu\left(\frac{V_z}{V}\right)_z.
			\label{2-1}
		\end{equation}
		The stress condition in \eqref{1-6} implies $V_z/V\to0$ as $z \to \pm \infty$.
		Then, integrating \eqref{2-1} over $\mathbb{R}$ gives
		\begin{equation}
			s^2=\frac{p(0)-p(v_-)}{v_-}.
			\label{2-2}
		\end{equation}
		Next, integrating \eqref{2-1} once from $z$ to $+\infty$ yields
		\begin{equation*}
			s^2V+p(V)+s\mu\frac{V_z}{V}=p(0).
		\end{equation*}
		It follows that
		\begin{equation}
			V_z=f(V):=-\frac{V}{s\mu}g(V),\qquad
			g(V)=p(V)-p(0)+s^2V.
			\label{2-4}
		\end{equation}
		Since $g''(V)=p''(V)>0$ and $g(0)=g(v_-)=0$, by strict convexity,
		\begin{equation}
			g(V)<0, \quad \hbox{for } 0<V<v_-.
			\label{2-5}
		\end{equation}
		Uniqueness for the scalar autonomous equation prevents a nonconstant
		connecting orbit from crossing either equilibrium. The prescribed end states then imply $V_{z} < 0$. Since $V>0$ and $\mu>0$, by \eqref{2-4}, we get
		\begin{equation}
			s<0.
			\label{2-6}
		\end{equation}
		Strict convexity of $p$ gives
		\begin{equation*}
			p'(0)<\frac{p(v_-)-p(0)}{v_-}<p'(v_-).
		\end{equation*}
		Using \eqref{2-2} and \eqref{2-6}, we obtain \eqref{1-8}.
		
		Conversely, suppose that \eqref{1-7} and \eqref{1-8} hold with $s<0$.
		Then, for $f$ defined in \eqref{2-4}, we have $f<0$ on $(0,v_-)$. Thus, the ODE \eqref{2-4} can be rewritten as
		\begin{align}
			\frac{\mathrm{d}z}{\mathrm{d}V} = \frac{1}{f(V)}.
			\label{2-81}
		\end{align}
		We choose any $v_*\in(0,v_-)$ and $z_*\in\R$ such that $V(z_*)=v_*$ and $z(v_*)=z_*$. In this case, we solve \eqref{2-81} in the form
		\begin{equation*}
			z-z_*=F(V):=\int_{v_*}^V\frac{\mathrm{d}y}{f(y)}.
		\end{equation*}
		Since $F'(V)=1/f(V)<0$, the function $F$ is strictly decreasing on $(0,v_-)$ and remains finite throughout this interval. Taylor's expansion of \eqref{2-4} at $V=0$ yields
		\begin{equation}
			f(V)=- \frac{V^{2}}{s\mu} \left[\frac{p(V)-p(0)}{V} + s^2\right] = -a_+V^2+O(V^3).
			\label{2-10}
		\end{equation}
		On the other hand, using \eqref{2-2}, we have
		\begin{equation}
			\begin{aligned}
				f(V) = - \frac{V}{s\mu} [p(V) - p(v_{-}) + s^2 (V - v_{-})]
				= & - \frac{V (V - v_{-})}{s\mu} \left[\frac{p(V)-p(v_{-})}{V - v_{-}} + s^2\right]
				\\
				= & \lambda_-(V-v_-)+O((V-v_-)^2), \quad \hbox{as } V\to v_-.
			\end{aligned}
			\label{2-11}
		\end{equation}
		Thus,
		\begin{equation*}
			\lim_{V \to 0^+} F(V) = +\infty,
			\quad
			\lim_{V \to v_-} F(V) = -\infty.
		\end{equation*}
		Then $V$ is a strictly decreasing solution of \eqref{2-4}, with
		\begin{equation*}
			V(-\infty) = v_-, \quad V(+\infty) = 0.
		\end{equation*}
		Thus $F$ is a bijection from $(0,v_{-})$ on to $\mathbb{R}$. According to the inverse function theorem, there exists a unique continuous inverse function $F^{-1}$ such that
		\begin{equation*}
			V = F^{-1}(z-z_\star), \quad z\in(-\infty,+\infty).
		\end{equation*}
		By integrating the first equation of \eqref{1-5} over $(z,+\infty)$, the shock wave profile $U(z)$ is uniquely determined. Moreover, by virtue of \eqref{2-4} and $U_z = - s V_z < 0$, \eqref{1-9} holds. Furthermore, \eqref{2-4} and $g(0) = g(v_-) = 0$ imply
		\begin{align*}
			\frac{V_z}{V} \to 0,
			\quad \text{as }\; z \to \pm \infty.
		\end{align*}
		Thus, all conditions in \eqref{1-6} hold.
		
		Finally, as $z \to +\infty$, \eqref{2-10} implies
		\begin{equation*}
			-\frac{V_{z}}{V^{2}}=\left(\frac1V\right)_z=a_++O(V)\to a_+,
		\end{equation*}
		which implies $V(z) \sim (a_{+}z)^{-1} \to 0$.
		As $z \to -\infty$, put $Y=v_- - V$. \eqref{2-11} gives
		\begin{equation*}
			Y_z=\lambda_-Y+O(Y^2).
		\end{equation*}
		Therefore, $Y$ decays exponentially as $z\to-\infty$ by comparison. It follows that $Y_z/Y-\lambda_-=(\log Y - \lambda_- z)_z=O(Y)$ is integrable. Integrating it over $[z,z_{0}]$ ($z_{0}\in\mathbb{R}^{+}$) and then exponentiating yields
		\begin{equation*}
			Y(z)e^{-\lambda_-z}=Y(z_{0})e^{-\lambda_-z_{0}} \exp\left\{-\int_{z}^{z_{0}} O(Y)\right\}.
		\end{equation*}
		Thus, $Y(z)e^{-\lambda_-z}$ tends to a finite positive constant
		$A_-$ as $z \to -\infty$. The proof of \eqref{1-911} is complete.
	\end{proof}
	
	\begin{proof}[Proof of Theorem \ref{theorem1.2}]
		Define
		\begin{equation}
			B=U-\mu\frac{V_z}{V}.
			\label{2-15}
		\end{equation}
		Using \eqref{1-5}, we obtain
		\begin{equation}
			\begin{cases}
				-sV_z-B_z=\mu\left(\dfrac{V_z}{V}\right)_z,\\
				-sB_z+p(V)_z=0.
			\end{cases}
			\label{2-17}
		\end{equation}
		Since $p'(V)<0$, $s<0$ and $V_z<0$, we conclude that $B_z<0$. Uniqueness up to translation follows from the uniqueness of $V$ and the definition \eqref{2-15}.
	\end{proof}
	
	We investigate the elementary properties of the shock wave solution $V(z)$.
	\begin{lemma}\label{lemma3.1}
		Let $V(z)$ be the shock wave obtained in Theorem \ref{theorem1.1}. It holds that
		\begin{equation}
			|V_z| \leq C V^2, \quad |V_{z z}| \leq CV|V_z| \leq C V^3, \quad |V_{zzz}| \leq C(|V| |V_{zz}| + V_z^{2}) \leq C V^4.
			\label{3-16}
		\end{equation}
		Moreover,
		\begin{equation}
			\left|\frac{V_z}{V^2}\right|\leq C,\quad
			\left|\left(\frac{V_z}{V^2}\right)_z\right|
			+\left|\left(\frac{V_z}{p'(V)V^2}\right)_z\right|
			\leq C|V_z|,
			\label{profile:coeff}
		\end{equation}
		for any $z \in(-\infty,+\infty)$.
	\end{lemma}
	\begin{proof}
		From \eqref{2-4} and the fact that $0<V \leq v_-$, there exists a constant $C > 0$ such that $|V_z|\leq CV^2$. Since $|f'(V)|\leq CV$ and $|f''(V)|\leq C$, we have
		\begin{equation*}
			|V_{zz}| = |f'(V)V_{z}| \leq CV|V_z| \leq C V^3, \quad |V_{zzz}| = |f'(V)V_{zz}+f''(V)V_{z}^{2}| \leq C(|V| |V_{zz}| + V_z^{2}) \leq C V^4.
		\end{equation*}
		For \eqref{profile:coeff}, based on \eqref{2-4}, we get
		\begin{equation}
			\frac{V_z}{V^{2}} = - \frac{1}{s\mu} \left(\frac{p(V) - p(0)}{V} + s^2\right) =-\frac{1}{s\mu}\left(\int_0^1p'(\theta V)\,\mathrm{d}\theta+s^2\right).
			\label{3.22}
		\end{equation}
		Hence $|V_{z}/V^{2}| \leq C$. By the regularity assumption $p \in C^{4}([0,\infty))$, it follows that
		\begin{equation}
			\begin{aligned}
				\left|\left(\frac{V_{z}}{V^{2}}\right)_{z}\right| = \left|\frac{V_{z}}{s\mu} \frac{d}{dV} \left(\int_0^1p'(\theta V)\,\mathrm{d}\theta+s^2\right)\right| = \left|\frac{V_{z}}{s\mu} \int_0^1 \theta p''(\theta V) \mathrm{d}\theta\right| \leq C |V_{z}|.
			\end{aligned}
			\label{3.23}
		\end{equation}
		Furthermore, combining \eqref{3-16} with \eqref{3.23}, we derive that
		\begin{equation*}
			\left|\left(\frac{V_{z}}{p'(V) V^{2}}\right)_{z}\right| = \left|\frac{1}{p'(V)} \left(\frac{V_{z}}{V^{2}}\right)_{z} - \frac{V_{z}}{V^{2}} \frac{p''(V) V_{z}}{(p'(V))^{2}}\right| \leq C |V_{z}|.
		\end{equation*}
		The proof is completed.
	\end{proof}

	\section{Nonlinear asymptotic stability}\label{sec:stability}
	\subsection{Perturbation equations and the main estimates}
	Let $(V,B)(x-st)$ be the viscous shock profile constructed in Theorem \ref{theorem1.2}. By the conservation of mass and \eqref{3-2}, we have
	\begin{equation}
		\int \binom{v(x,t)-V(x+x_{0}-st)}{b(x,t)-B(x+x_{0}-st)} \mathrm{d} x
		= \int \binom{v_{0}(x)-V_{0}(x)}{b_{0}(x)-B_{0}(x)} \mathrm{d} x
		= \binom{0}{0},
		\label{3-21}
	\end{equation}
	which makes it possible to introduce the anti--derivative technique to study the asymptotic stability of the shock wave. We define
	\begin{equation}
		(\phi, \psi)(z, t)=\int_{-\infty}^z\left(v(y + st, t)-V\left(y+x_0\right), b(y + st, t)-B\left(y+x_0\right)\right) \mathrm{d} y,
		\label{3-3}
	\end{equation}
	where $z=x-s t$. By \eqref{3-21}, we obtain
	\begin{equation}
		(\phi, \psi)( \pm \infty, t)=0 \quad \hbox{for all } t \geq 0.
		\label{3-4}
	\end{equation}
	Then, \eqref{3-3} is equivalent to
	\begin{align}
		(v, b)(x, t)=(V, B)\left(x+x_0-s t\right)+\left(\phi_z, \psi_z\right)(z, t).
		\label{3-5}
	\end{align}
	Substituting \eqref{3-5} into \eqref{1-11}, making use of \eqref{2-17} and \eqref{3-4}, and integrating the resulting equations over $(-\infty,z)$, we obtain the following perturbation equations
	\begin{equation}
		\begin{cases}
			\phi_t=s \phi_z+\psi_z+\mu\left(\frac{V_z+\phi_{z z}}{V+\phi_z} - \frac{V_z}{V}\right),
			\\
			\psi_t=s \psi_z - (p(V+\phi_z)-p(V)).
		\end{cases}
		\label{3-611}
	\end{equation}
	Without loss of generality, we normalize the translation by taking $x_0=0$. The initial data are
	\begin{align}
		(\phi_0, \psi_0)(z) := (\phi, \psi)(z, 0)=\int_{-\infty}^z (v_0-V, b_0-B)(y) \mathrm{d} y,
		\label{3-9}
	\end{align}
	with $(\phi_0, \psi_0)( \pm \infty)=0$ due to \eqref{3-2}.
	
	We define
	\begin{equation}
		r=\phi_{z}, \quad q=\phi_t-s\phi_z=\psi_z+\mu\left(\frac{V_z+\phi_{z z}}{V+\phi_z} - \frac{V_z}{V}\right).
		\label{defs:r-q}
	\end{equation}
	By \eqref{1-10} and \eqref{2-15}, $q$ is exactly the perturbation of the original velocity. The initial values $r_{0}$ and $q_{0}$ are prescribed by
	\begin{equation}
		r_{0}=\phi_{0z}, \quad q_0=\psi_{0z} + \mu\left(\frac{V_z+\phi_{0z z}}{V+\phi_{0z}} - \frac{V_z}{V}\right).
		\label{q:init}
	\end{equation}
	The weight functions $w_{i}$ are specified as follows
	\begin{equation}
		w_i=V^{-i},\quad i=1,2,3.
		\label{3-10}
	\end{equation}
	For any fixed $T>0$, we define the solution space
	\[
	X(0,T)=\{(\phi,\psi):\
	\phi\in C([0,T];L^2),\ \psi\in C([0,T];H^2),\
	r,q\in C([0,T];H^1_{w_2})\}.
	\]
	Denote
	\begin{equation}
		\mathcal E(t)=\|\phi(t)\|^2+\|\psi(t)\|_2^2
		+\|r(t)\|_{1,w_2}^2+\|q(t)\|_{1,w_2}^2,\quad
		N(t)=\sup_{0\leq\tau\leq t}\mathcal E(\tau)^{1/2}.
		\label{3-11}
	\end{equation}
	
	\begin{proposition}[Local existence]\label{proposition3.1}
		Suppose that $\phi_0\in L^2(\mathbb{R})$, $\psi_0\in H^2(\mathbb{R})$ and
		$r_{0},q_0\in H^1_{w_2}(\mathbb{R})$, with $r_{0}$ and $q_0$ given by \eqref{q:init}. For any sufficiently small $\delta_0 > 0$, there exists a positive constant $T_0$ depending on $\delta_0$ such that if $N(0) \leq \delta_0$, then the problem \eqref{3-611} with \eqref{3-9} has a unique solution $(\phi, \psi) \in X(0, T_0)$ satisfying $N(t) \leq C_{*} N(0)$ for any $0 \leq t \leq T_0$, where $C_{*}$ is a constant independent of $t$.
	\end{proposition}
	The proof is given in Appendix~\ref{sec:local}. We define the positive coefficient
	\[
	d(V)=\frac{s p''(V)V_z}{p'(V)^2}.
	\]
	Since $p \in C^{4}([0,\infty))$, we get $c|V_z|\leq d(V)\leq C|V_z|$.
	The weighed energy estimates used below are
	\begin{equation*}
		\mathcal D(t) =\|r(t)\|_{2,w_1}^2+\|q(t)\|_{w_1}^2
		+\|q_z(t)\|_{1,w_3}^2+
		\int d(V)(\psi^2+\psi_z^2+\psi_{zz}^2)(t).
	\end{equation*}
	
	\begin{proposition}[A priori estimate]\label{proposition3.2}
		Let the assumptions of Proposition \ref{proposition3.1} hold. Assume that $(\phi, \psi) \in X(0, T)$ is a solution of \eqref{3-611} with \eqref{3-9} for a positive constant $T$. Then there exists a constant $\delta_{1}$ satisfying $0 < \delta_{1} < \delta_{0}$, independent of $T$, such that if
		\begin{equation}
			N(t) \leq \delta_1,
			\label{3-14}
		\end{equation}
		for any $0 \leq t \leq T$, it holds that
		\begin{equation}
			\mathcal E(t)+\int_0^t\mathcal D(\tau)\,\mathrm{d}\tau
			\leq C\mathcal E(0).
			\label{3-12}
		\end{equation}
	\end{proposition}
	The Sobolev embedding inequality $\|f\|_{L^{\infty}} \leq \sqrt{2} \|f\|^{1/2}\|f_{z}\|^{1/2}$ and \eqref{3-11} yield
	\begin{equation}
		\sup _{\tau \in[0, t]}\left\{\|(\phi, \psi, \psi_z)(\cdot, \tau)\|_{L^{\infty}}, \left\|(r/V, q/V)(\cdot, \tau)\right\|_{L^{\infty}}\right\} \leq C N(t),
		\label{3-15}
	\end{equation}
	for some constant $C$ independent of $\delta_1$ and $T$. In addition, by the definition of $q$ in \eqref{defs:r-q}, we have
	\begin{equation}
		\frac{r_{z}}{V}= \frac{r+V}{V} \frac{q-\psi_z}{\mu} + \frac{V_z(V+r)}{V^2} - \frac{V_z}{V} = \left(1+ \frac{r}{V}\right) \frac{q-\psi_z}{\mu} + \frac{V_z r}{V^2}.
		\label{3-161}
	\end{equation}
	By virtue of \eqref{3-16}, it follows that
	\begin{equation}
		\|r_{z}/V\|_{L^\infty}\leq CN(t).
		\label{rz:infinity}
	\end{equation}
	Using \eqref{3-15}, when $N(t)$ is small enough, it is obvious that
	\begin{equation}
		|r/V| \leq 1/2, \quad 1 + r/V > 0, \quad V + r > 0.
		\label{3-151}
	\end{equation}
	We define
	\begin{equation}
		\eta = \log(1+r/V),
		\label{3-6111}
	\end{equation}
	which implies
	\begin{equation}
		\eta_{z} = (\log(1+r/V))_{z} = \frac{V_z+r_{z}}{V+r} - \frac{V_z}{V}.
		\label{3-6112}
	\end{equation}
	It follows that
	\begin{equation}
		q=\phi_t-s\phi_z=\psi_z+\mu\eta_z.
		\label{3-61}
	\end{equation}
	System \eqref{3-611} can therefore be written as
	\begin{equation}
		\begin{cases}
			\phi_t=sr+\psi_z+\mu(r/V)_z+\mu W_z,
			\\
			\psi_t=s\psi_z-p'(V)r-G,
		\end{cases}
		\label{3-6}
	\end{equation}
	where
	\begin{equation}
		W = \log(1+r/V)-r/V,
		\label{3-7}
	\end{equation}
	and
	\begin{equation}
		G=p(V+r)-p(V)-p'(V)r.
		\label{3-8}
	\end{equation}
	The local existence and the a priori estimate yield the following global stability theorem.
	
	\begin{theorem}\label{theorem3.1}
		Under the assumptions of Proposition \ref{proposition3.1}, there exists
		a constant $\delta_2>0$ such that if $N(0)\leq\delta_2$, the Cauchy problem \eqref{3-611} with \eqref{3-9} has a unique global solution $(\phi, \psi) \in X(0, \infty)$ satisfying \eqref{3-12} for any $t \in [0, \infty)$. Moreover, it follows that
		\begin{equation}
			\|r(t)\|_{L^\infty}+\|q(t)\|_{L^\infty} \to 0,
			\quad\text{as }t\to\infty.
			\label{3-13}
		\end{equation}
	\end{theorem}

	\subsection{Basic and interactive energy estimates}
	\begin{lemma}\label{lemma3.2}
		Under the a priori assumption \eqref{3-14}, if $N(T)$ is sufficiently small, then it holds that
		\begin{equation}
			\begin{aligned}
				|W|&\leq Cr^2/V^2,& |G|&\leq Cr^2,\\
				|W_z|&\leq C\left(\frac{|V_z|}{V^3} r^2+\frac{|r||r_{z}|}{V^2}\right),&
				|G_z|&\leq C(|V_z|r^2+|r r_z|),\\
				|W_{zz}|&\leq C\left[\left(\frac{V_z^2}{V^4} + \frac{|V_{zz}|}{V^3}\right)r^2 + \frac{r_{z}^2}{V^2} + \frac{|V_z|}{V^3}|r||r_{z}| + \frac{|r||r_{zz}|}{V^2}\right],\\
				|G_{zz}|&\leq C (V_z^2r^2 + |V_z||r||r_{z}| + |V_{zz}|r^2 + r_{z}^2 + |r||r_{zz}|).
			\end{aligned}
			\label{3-18}
		\end{equation}
	\end{lemma}
	\begin{proof}
		From \eqref{3-7} and \eqref{3-8}, a direct calculation gives $|W| \leq C r^{2}/V^{2}$ and $|G| \leq C r^2$. Furthermore, it holds that
		\begin{align*}
			W_z =\frac{V_z+r_{z}}{V+r}-\frac{V_z}{V}-\frac{r_{z}}{V}+\frac{V_z}{V^2} r = \left(\frac{1}{V+r}-\frac{1}{V}+\frac{r}{V^2}\right) V_z + \left(\frac{1}{V+r}-\frac{1}{V}\right) r_{z},
		\end{align*}
		and
		\begin{align*}
			G_z & =p'(V+r)(V_z+r_{z})-p'(V) V_z-p''(V) V_z r-p'(V) r_{z}
			\\
			& =(p'(V+r)-p'(V)-p''(V) r) V_z+(p'(V+r)-p'(V)) r_{z}.
		\end{align*}
		Taylor's expansion then yields the estimates of $W_z$ and $G_z$. An analogous argument gives
		\begin{align*}
			W_{zz} = & \left(\frac{V_z+r_{z}}{V+r} - \frac{V_z}{V} - \frac{r_{z}}{V} + \frac{V_z}{V^2}r\right)_z
			\\
			= & \frac{V_{zz}+r_{zz}}{V+r} - \frac{(V_z+r_{z})^2}{(V+r)^2} + \frac{V_z^2}{V^2} - \frac{V_{zz}}{V} + \frac{2V_z}{V^2}r_{z} - \frac{r_{zz}}{V} - \frac{2V_z^2}{V^3}r + \frac{V_{zz}}{V^2}r
			\\
			= & \left[-\frac{1}{(V+r)^2} - \left(-\frac{1}{V^2}\right) -\frac{2r}{V^3}\right]V_z^2 + \left[\frac{1}{V+r} - \frac{1}{V} - \left(-\frac{r}{V^2}\right)\right]V_{zz}
			\\
			& - \frac{r_{z}^2}{(V+r)^2} + 2\left[-\frac{1}{(V+r)^2} - \left(-\frac{1}{V^2}\right)\right]V_z r_{z} + \left(\frac{1}{V+r} - \frac{1}{V}\right)r_{zz},
		\end{align*}
		and
		\begin{align*}
			G_{zz} = & \left[\left(p'(V+r)-p'(V)-p''(V)r\right)V_z + \left(p'(V+r)-p'(V)\right)r_{z}\right]_z
			\\
			= & \left[p''(V+r)(V_z+r_{z}) - p''(V)V_z - p'''(V)V_z r - p''(V)r_{z}\right]V_z
			\\
			& + \left(p'(V+r)-p'(V)-p''(V)r\right)V_{zz} + \left[
			p''(V+r)(V_z+r_{z})
			\right.
			\\
			& \left. \quad -p''(V)V_z\right]r_{z} + \left(p'(V+r)-p'(V)\right)r_{zz}
			\\
			= & \left(p''(V+r)-p''(V)-p'''(V)r\right)V_z^2 + \left(p''(V+r)-p''(V)\right)r_{z}V_z
			\\
			& + \left(p'(V+r)-p'(V)-p''(V)r\right)V_{zz} + \left(p''(V+r)-p''(V)\right)V_z r_{z}
			\\
			& + p''(V+r)r_{z}^2 + \left(p'(V+r)-p'(V)\right)r_{zz}.
		\end{align*}
		Using \eqref{3-15}, $p \in C^{4}$ and Taylor's expansion, we can get the estimates of $W_{zz}$ and $G_{zz}$.
	\end{proof}
	
	\begin{lemma}\label{lemma3.3}
		Let the assumptions of Theorem \ref{theorem3.1} hold. If $N(T) \ll 1$, then there exists a constant $C>0$ independent of $T$ such that
		\begin{equation}
			\|\phi\|^2+\|\psi\|^2+\int_0^t \|r\|_{w_1}^2 + \int_0^t \int d(V) \psi^2 \leq C (\|\phi_0\|^2 + \|\psi_0\|^2),
			\label{3-22}
		\end{equation}
		for any $t \in[0, T]$.
	\end{lemma}
	\begin{proof}
		Multiplying \eqref{3-6}$_{1}$ by $\phi$ and \eqref{3-6}$_{2}$ by $\frac{\psi}{-p'(V)}$, respectively, adding the resulting equations and integrating over $\mathbb{R} \times [0, t]$ yields
		\begin{equation}
			\begin{aligned}
				& \frac{1}{2} \int \left(\phi^2-\frac{\psi^2}{p'(V)}\right) + \mu \int_0^t \int \frac{r^2}{V} + \frac{1}{2} \int_0^t \int d(V) \psi^2
				\\
				= & \frac{1}{2} \int \left(\phi_0^2-\frac{\psi_0^2}{p'(V)}\right) + \int_0^t \int \left(\mu \frac{\phi r}{V} + \frac{s}{2} \phi^2 - \frac{s \psi^2}{2 p'(V)} + \phi \psi + \mu W \phi\right)_z
				\\
				& - \mu \int_0^t \int W r + \int_0^t \int G \frac{\psi}{p'(V)}.
			\end{aligned}
			\label{3-23}
		\end{equation}
		By \eqref{3-11}, the functions $\phi,\psi,r/V$ belong to $H^1$ and hence tend to zero as $z \to \pm \infty$. Since $|W| \leq C r^{2}/V^{2}$, we get $W \to 0$ as $z \to \pm \infty$. It then holds that
		\begin{align}
			\int_0^t \int \left(\mu \frac{\phi r}{V} + \frac{s}{2} \phi^2 - \frac{s \psi^2}{2 p'(V)} + \phi \psi + \mu W \phi\right)_z = 0.
			\label{3-231}
		\end{align}
		By virtue of \eqref{3-15} and \eqref{3-18}, we obtain
		\begin{equation}
			\begin{aligned}
				\left| -\mu \int_0^t \int W r + \int_0^t \int G \frac{\psi}{p'(V)} \right| \leq C \int_0^t \int \frac{|r|^3}{V^2} + C \int_0^t \int r^2 |\psi| \leq C N(t) \int_0^t \int \frac{r^2}{V},
			\end{aligned}
			\label{3-24}
		\end{equation}
		where we have used the elementary bound $1 \leq C/V$ for $0 < V \leq v_-$. Inserting \eqref{3-231} and \eqref{3-24} into \eqref{3-23}, we deduce that
		\begin{align*}
			\frac{1}{2} \int \left(\phi^2-\frac{\psi^2}{p'(V)}\right) + (\mu - C N(t)) \int_0^t \int \frac{r^2}{V} + \frac{1}{2} \int_0^t \int d(V) \psi^2 \leq \frac{1}{2} \int \left(\phi_0^2-\frac{\psi_0^2}{p'(V)}\right),
		\end{align*}
		which implies \eqref{3-22}.
	\end{proof}
	
	\begin{lemma}[Energy estimates for velocity]\label{lemma3.4}
		Let the assumptions of Theorem \ref{theorem3.1} hold. If $N(T) \ll 1$, then there exists a constant $C>0$ independent of $T$ such that
		\begin{equation}
			\| r \|_{w_2}^2 + \int_0^t \|q\|_{w_1}^2 \leq C (\|\phi_0\|^2 + \|\psi_0\|^2 + \|r_{0}\|_{w_2}^2),
			\label{3-25}
		\end{equation}
		for any $t \in[0, T]$.
	\end{lemma}
	\begin{proof}
		To obtain the energy estimates of $q$, we exploit the coupling between the two perturbation equations. By \eqref{defs:r-q}, it is obvious that
		\begin{equation}
			r_{t} = r_{t} - s r_{z} + s r_{z} = q_{z} + s r_{z}.
			\label{3-251}
		\end{equation}
		Multiplying \eqref{3-251} by $\psi/V$ and \eqref{3-611}$_{2}$ by $r/V$, respectively, adding the resulting equations and integrating over $\mathbb{R}$ yields
		\begin{equation}
			\begin{aligned}
				\frac{\mathrm{d}}{\mathrm{d}t} \int \frac{r}{V} \psi = \int \left(\frac{r_{t}}{V} \psi + \frac{r}{V} \psi_{t}\right) = & \int \frac{q_{z}}{V} \psi + s \int \left(\frac{r}{V}\right)_{z} \psi + s \int \frac{V_{z}}{V^{2}} r \psi
				\\
				& + s \int \frac{r}{V} \psi_{z} - \int \frac{r}{V} [p(V + r) - p(V)].
			\end{aligned}
			\label{3.3}
		\end{equation}
		Since $\psi,r/V,q/V \in H^{1}$, we have
		\begin{equation}
			\psi,r/V,q/V \to 0, \quad \hbox{as } z \to \pm \infty,
			\label{3.1211}
		\end{equation}
		which implies
		\begin{equation}
			\eta \to 0, \quad \hbox{as } z \to \pm \infty.
			\label{3.121}
		\end{equation}
		By \eqref{3.1211}, we know
		\begin{equation}
			s \int \left(\frac{r}{V}\right)_{z} \psi + s \int \frac{r}{V} \psi_{z} = s \int \left(\frac{r}{V} \psi\right)_{z} = 0,
			\label{3.311}
		\end{equation}
		and
		\begin{equation}
			\begin{aligned}
				\int \frac{q_{z}}{V} \psi = \int \left(\frac{q \psi}{V}\right)_{z} - \int \frac{q}{V} \psi_{z} + \int \frac{V_{z}}{V^{2}} q \psi = - \int \frac{q^{2}}{V} + \mu \int \frac{q}{V} \eta_{z} + \int \frac{V_{z}}{V^{2}} q \psi,
			\end{aligned}
			\label{3.3111}
		\end{equation}
		where we have used $\psi_z = q - \mu \eta_{z}$ in the last equality. Combining \eqref{3.311} with \eqref{3.3111}, \eqref{3.3} becomes
		\begin{equation}
			\begin{aligned}
				\frac{\mathrm{d}}{\mathrm{d}t} \int \frac{r}{V} \psi = - \int \frac{q^{2}}{V} + s \int \frac{V_{z}}{V^{2}} r \psi + \mu \int \frac{\eta_{z} q}{V} + \int \frac{V_{z}}{V^{2}} q \psi - \int \frac{r}{V} [p(V + r) - p(V)].
			\end{aligned}
			\label{3.4}
		\end{equation}
		However, the estimates in Lemma \ref{lemma3.3} do not control the term $\mu \int \eta_{z} q/V$. Therefore, we introduce an auxiliary functional $\Phi(\eta)$ to eliminate this term. We define
		\begin{align*}
			\Phi(\eta) = \mu [e^{\eta} (\eta - 1) + 1].
		\end{align*}
		Then
		\begin{align}
			\Phi'(\eta) = \mu e^{\eta} (\eta - 1) + \mu e^{\eta} = \mu \eta e^{\eta}.
			\label{3.6}
		\end{align}
		As $r = V (e^{\eta} - 1)$, we have
		\begin{equation*}
			r_{t}=q_{z}+sr_{z}=q_{z}+sV_{z}(e^{\eta}-1)+sVe^{\eta}\eta_{z},
		\end{equation*}
		and thus,
		\begin{align*}
			Ve^{\eta}\eta_{t}=q_{z}+sV_{z}(e^{\eta}-1)+sVe^{\eta}\eta_{z}.
		\end{align*}
		By virtue of \eqref{3.6}, we obtain
		\begin{equation}
			\frac{\mathrm{d}}{\mathrm{d}t} \int \Phi(\eta) = \int \Phi'(\eta) \eta_{t} = \mu \int \eta e^{\eta} \eta_{t} = s \mu \int \eta e^{\eta} \eta_{z} + s \mu \int \frac{V_{z}}{V} \eta (e^{\eta} - 1) + \mu \int \frac{\eta q_{z}}{V}.
			\label{3.12}
		\end{equation}
		Integration by parts, together with \eqref{3.1211} and \eqref{3.121}, gives
		\begin{align*}
			\mu \int \frac{\eta q_{z}}{V} = \mu \int \left(\frac{\eta q}{V}\right)_{z} - \mu \int \frac{\eta_{z} q}{V} + \mu \int \frac{V_{z}}{V^{2}} \eta q = - \mu \int \frac{\eta_{z} q}{V} + \mu \int \frac{V_{z}}{V^{2}} \eta q.
		\end{align*}
		Moreover, since $\Phi(\eta) \to 0$ as $z \to \pm \infty$, utilizing $\Phi_z(\eta)=\mu\eta e^{\eta} \eta_{z}$, we derive
		\begin{align*}
			s \mu \int \eta e^{\eta} \eta_{z} = s \int \Phi_{z}(\eta) = 0.
		\end{align*}
		Using $e^{\eta}-1 = r/V$, we rewrite \eqref{3.12} as
		\begin{equation}
			\begin{aligned}
				\frac{\mathrm{d}}{\mathrm{d}t} \int \Phi(\eta) = & s \mu \int \frac{V_{z}}{V} \eta (e^{\eta} - 1) - \mu \int \frac{\eta_{z} q}{V} + \mu \int \frac{V_{z}}{V^{2}} \eta q
				\\
				= & s \mu \int \frac{V_{z}}{V^{2}} \eta r - \mu \int \frac{\eta_{z} q}{V} + \mu \int \frac{V_{z}}{V^{2}} \eta q.
			\end{aligned}
			\label{3.15}
		\end{equation}
		Adding \eqref{3.4} and \eqref{3.15}, we deduce that
		\begin{align*}
			& \frac{\mathrm{d}}{\mathrm{d}t} \int \Phi(\eta) + \frac{\mathrm{d}}{\mathrm{d}t} \int \frac{r}{V} \psi + \int \frac{q^{2}}{V}
			\\
			= & s \int \frac{V_{z}}{V^{2}} r \psi + \int \frac{V_{z}}{V^{2}} q (\psi + \mu \eta) - \int \frac{r}{V} [p(V + r) - p(V)] + s \mu \int \frac{V_{z}}{V^{2}} r \eta.
		\end{align*}
		Since $q = (\psi + \mu \eta)_{z}$, integration by parts gives
		\begin{equation*}
			\int \frac{V_{z}}{V^{2}} q (\psi + \mu \eta) = \frac{1}{2} \int \left[\frac{V_{z}}{V^{2}} (\psi + \mu \eta)^{2}\right]_{z} - \frac{1}{2} \int \left(\frac{V_{z}}{V^{2}}\right)_{z} (\psi + \mu \eta)^{2} = - \frac{1}{2} \int \left(\frac{V_{z}}{V^{2}}\right)_{z} (\psi + \mu \eta)^{2}.
		\end{equation*}
		Combining the above identity, we conclude that
		\begin{equation}
			\begin{aligned}
				\frac{\mathrm{d}}{\mathrm{d}t} \int \Phi(\eta) + \frac{\mathrm{d}}{\mathrm{d}t} \int \frac{r}{V} \psi + \int \frac{q^{2}}{V} = & s \int \frac{V_{z}}{V^{2}} r \psi + s \mu \int \frac{V_{z}}{V^{2}} r \eta - \frac{1}{2} \int \left(\frac{V_{z}}{V^{2}}\right)_{z} (\psi + \mu \eta)^{2}
				\\
				& - \int \frac{r}{V} [p(V + r) - p(V)].
			\end{aligned}
			\label{3.18}
		\end{equation}
		Next, we need to eliminate the unfavorable term $s \int r \psi V_{z}/V^{2}$. Multiplying \eqref{3-6}$_{2}$ by $\frac{sV_{z}}{p'(V)V^{2}}\psi$ yields
		\begin{equation}
			\frac{\mathrm{d}}{\mathrm{d}t} \int \frac{s V_{z}}{2 p'(V) V^{2}} \psi^{2} = \int \frac{s^{2} V_{z}}{p'(V) V^{2}} \psi \psi_{z} - s \int \frac{V_{z}}{V^{2}} \psi r - s \int \frac{V_{z}}{p'(V) V^{2}} \psi G,
			\label{3.20}
		\end{equation}
		where, noting \eqref{3.1211} and using integration by parts, we obtain
		\begin{equation}
			\begin{aligned}
				\int \frac{s^{2} V_{z}}{p'(V) V^{2}} \psi \psi_{z} = \int \left(\frac{s^{2} V_{z} \psi^{2}}{2 p'(V) V^{2}}\right)_{z} - \int \left(\frac{s^{2} V_{z}}{2 p'(V) V^{2}}\right)_{z} \psi^{2} = - \int \left(\frac{s^{2} V_{z}}{2 p'(V) V^{2}}\right)_{z} \psi^{2}.
			\end{aligned}
			\label{3.201}
		\end{equation}
		Combining \eqref{3.18}, \eqref{3.20} and \eqref{3.201}, we get
		\begin{equation}
			\begin{aligned}
				& \frac{\mathrm{d}}{\mathrm{d}t} \int \Phi(\eta) + \frac{\mathrm{d}}{\mathrm{d}t} \int \frac{r}{V} \psi + \frac{\mathrm{d}}{\mathrm{d}t} \int \frac{s V_{z}}{2 p'(V) V^{2}} \psi^{2} + \int \frac{q^{2}}{V}
				\\
				= & s \mu \int \frac{V_{z}}{V^{2}} r \eta - \frac{1}{2} \int \left(\frac{V_{z}}{V^{2}}\right)_{z} (\psi + \mu \eta)^{2} - \int \frac{r}{V} [p(V + r) - p(V)]
				\\
				& - \int \left(\frac{s^{2} V_{z}}{2 p'(V) V^{2}}\right)_{z} \psi^{2} + s \int \frac{V_{z}}{-p'(V) V^{2}} \psi G
				\\
				=: & \sum_{i=1}^{5} I_{i}.
			\end{aligned}
			\label{3.21}
		\end{equation}
		Now we estimate the five terms $\{I_{i}\}_{i=1}^{5}$ on the right--hand side of \eqref{3.21}. By Taylor's expansion, it holds that
		\begin{equation}
			\left| \eta \right| = \left| \log \left(1 + \frac{r}{V}\right) \right| \leq C \left|\frac{r}{V}\right|.
			\label{3.242}
		\end{equation}
		Utilizing \eqref{3-16} and \eqref{3.242}, we have
		\begin{equation}
			|I_{1}| = \left| s \mu \int \frac{V_{z}}{V^{2}} r \eta \right| \leq C \int \frac{r^{2}}{V}.
			\label{3.25}
		\end{equation}
		For $I_{2} + I_{4}$, the estimates \eqref{profile:coeff} and \eqref{3.242} lead to
		\begin{align*}
			|I_{2} + I_{4}| \leq C \int |V_{z}| [(\psi + \mu \eta)^{2} + \psi^{2}] \leq C \int r^{2} + C \int |V_{z}| \psi^{2}.
		\end{align*}
		Using \eqref{3-16}, \eqref{3-18}, Taylor's expansion and the fact $\| \psi \|_{L^{\infty}} \leq N(t)$, we estimate \(I_3+I_5\) as
		\begin{equation}
			|I_{3} + I_{5}| \leq C \int \frac{r^{2}}{V} + C N(t) \int r^{2}.
			\label{3.27}
		\end{equation}
		Substituting \eqref{3.25}--\eqref{3.27} into \eqref{3.21} and organizing the terms, we infer that
		\begin{equation}
			\frac{\mathrm{d}}{\mathrm{d}t} \int \Phi(\eta) + \frac{\mathrm{d}}{\mathrm{d}t} \int \frac{r}{V} \psi + \frac{\mathrm{d}}{\mathrm{d}t} \int \frac{s V_{z}}{2 p'(V) V^{2}} \psi^{2} + \int \frac{q^{2}}{V} \leq C \int \frac{r^{2}}{V} + C \int r^{2} + C \int |V_{z}| \psi^{2}.
			\label{3.28}
		\end{equation}
		Integrating \eqref{3.28} over $[0,t]$ and using $1 \leq C V^{-1}$ for $0<V\leq v_-$, we get
		\begin{align*}
			& \int \Phi(\eta) + \int \frac{r}{V} \psi + \int \frac{s V_{z}}{2 p'(V) V^{2}} \psi^{2} + \int_{0}^{t} \int \frac{q^{2}}{V}
			\\
			\leq & \int \left(\Phi(\eta) + \frac{r}{V} \psi\right)(z, 0) + \int \frac{s V_{z}}{2 p'(V) V^{2}} \psi_{0}^{2} + C \int_{0}^{t} \int \frac{r^{2}}{V} + C \int_{0}^{t} \int |V_{z}| \psi^{2}.
		\end{align*}
		Since $s < 0$, $V_{z} < 0$ and $p'(V) < 0$, it follows that $\int \frac{s V_{z}}{2 p'(V) V^{2}} \psi^{2} < 0$. By $p \in C^{4}([0,\infty))$, we have
		\begin{equation}
			\left|\int \frac{s V_{z}}{2 p'(V) V^{2}} \psi^{2}\right| \leq C \int \psi^{2}.
			\label{3.291}
		\end{equation}
		Define $K(\theta) = \Phi(\log(1+\theta))$, where $\theta = r/V$. Differentiating it with respect to $\theta$ gives
		\begin{equation*}
			K'(\theta) = \mu \log (1 + \theta), \quad K''(\theta) = \mu (1 + \theta)^{-1}.
		\end{equation*}
		By virtue of $|\theta| \leq |r/V| \leq 1/2$ in \eqref{3-151}, we see
		\begin{equation*}
			\frac{2}{3} \mu \leq K''(\theta) \leq 2 \mu.
		\end{equation*}
		Since $K(0) = K'(0) = 0$, Taylor's theorem gives
		\begin{equation*}
			K(\theta) = \frac{1}{2} K''(\xi) \cdot \theta^{2} = \frac{1}{2} K''(\xi) \cdot \left(\frac{r}{V}\right)^{2}, \quad \hbox{for some } \xi \in (0,\theta).
		\end{equation*}
		Consequently,
		\begin{equation}
			\frac{1}{3} \mu \left(\frac{r}{V}\right)^{2} \leq \Phi(\eta) \leq \mu \left(\frac{r}{V}\right)^{2}.
			\label{3.34}
		\end{equation}
		Applying Young's inequality gives
		\begin{equation}
			\left| \int \frac{r}{V} \psi \right| \leq \frac{\mu}{6} \int \left(\frac{r}{V}\right)^{2} + C \int \psi^{2}.
			\label{3.35}
		\end{equation}
		Combining \eqref{3.291}, \eqref{3.34} and \eqref{3.35} with Lemma \ref{lemma3.3}, we derive that
		\begin{align*}
			& \int \left(\frac{r}{V}\right)^{2} + \int_{0}^{t} \int \frac{q^{2}}{V}
			\\
			\leq & C (\|\phi_0\|^2 + \|\psi_0\|^2 + \|r_{0}\|_{w_2}^2) + C \int \psi^{2} + C \int_{0}^{t} \int \frac{r^{2}}{V} + C \int_{0}^{t} \int |V_{z}| \psi^{2}
			\\
			\leq & C (\|\phi_0\|^2 + \|\psi_0\|^2 + \|r_{0}\|_{w_2}^2).
		\end{align*}
		The proof is completed.
	\end{proof}

	\subsection{Weighted velocity estimates}
	By \eqref{3-61}, we obtain
	\[
	q_{t}=\psi_{zt}+\mu\eta_{zt},\qquad sq_{z}=s\psi_{zz}+s\mu\eta_{zz}.
	\]
	Consequently,
	\begin{equation}
		q_t = sq_z + (\psi_{t} - s\psi_{z})_{z} + \mu (\eta_{t} - s\eta_{z})_{z}.
		\label{4.71}
	\end{equation}
	\eqref{3-611} implies $(\psi_{t} - s\psi_{z})_{z} = - [p(V+r)-p(V)]_{z}$. From \eqref{3-6111}, \eqref{3-6112} and $r_t=s r_z+q_z$,
	\begin{equation*}
		\begin{aligned}
			\eta_t-s\eta_z=\frac{r_{t}}{V+r}-s\left(\frac{V_z+r_{z}}{V+r}-\frac{V_z}{V}\right)=\frac{q_z}{V+r}
			+\frac{sV_zr}{V(V+r)}.
		\end{aligned}
	\end{equation*}
	We then rewrite \eqref{4.71} as
	\begin{equation}
		q_t = sq_z + \left[\mu\left(\frac{q_z}{V+r} + \frac{sV_zr}{V(V+r)}\right) - (p(V+r)-p(V))\right]_z.
		\label{4.7}
	\end{equation}
	As $(V+r)^{-1} = V^{-1} V/(V+r) = V^{-1} (1+r/V)^{-1}$, we get from \eqref{3-151} that
	\begin{equation}
		\frac{2}{3V} \leq \frac{1}{V+r} \leq \frac{2}{V}.
		\label{4.711}
	\end{equation}
	To handle the singularity of the equation at $z = +\infty$, we choose a cut--off function $\kappa \in C_c^\infty(\R)$ satisfying
	\begin{equation}
		\begin{cases}
			\kappa(z) = 1, \, \;\hbox{for } |z| < L; \quad \kappa(z) = 0, \, \;\hbox{for } |z| > 2L;
			\\
			0 \leq \kappa(z) \leq 1, \quad |\kappa_{z}(z)| \leq \frac{C}{L}, \,\; \hbox{for } z \in \mathbb{R}.
		\end{cases}
		\label{4.10}
	\end{equation}
	
	\begin{lemma}\label{lemma3.6}
		Let the assumptions of Theorem \ref{theorem3.1} hold. If $N(T) \ll 1$, then there exists a constant $C>0$ independent of $T$ such that
		\begin{equation}
			\|q(t)\|_{w_2}^2+\int_0^t\|q_z\|_{w_3}^2\,\mathrm{d}\tau
			\leq C\bigl(\|\phi_0\|^2+\|\psi_0\|^2
			+\|r_0\|_{w_2}^2+\|q_0\|_{w_2}^2\bigr),
			\label{3-46}
		\end{equation}
		for any $t \in[0, T]$.
	\end{lemma}
	\begin{proof}
		Multiplying both sides of \eqref{4.7} by $\frac{q}{V^{2}}$, we obtain
		\begin{equation}
			\begin{aligned}
				\frac{1}{2} \frac{\mathrm{d}}{\mathrm{d}t} \frac{q^{2}}{V^{2}} = & \left[\mu \left(\frac{q_z}{V+r} + \frac{sV_zr}{V(V+r)}\right)\frac{q}{V^{2}}  - (p(V+r)-p(V)) \frac{q}{V^{2}} + \frac{s}{2} \frac{q^{2}}{V^{2}}\right]_{z} + s \frac{V_{z} q^{2}}{V^{3}}
				\\
				& - \mu \frac{q_{z}}{V+r} \left(\frac{q_{z}}{V^{2}} - \frac{2 V_{z} q}{V^{3}}\right) - \bigg[\frac{s \mu V_zr}{V(V+r)} - (p(V+r)-p(V))\bigg] \left(\frac{q_{z}}{V^{2}} - \frac{2 V_{z} q}{V^{3}}\right).
			\end{aligned}
			\label{3-4611}
		\end{equation}
		Next, we multiply the preceding identity by $\kappa^{2}(z)$ defined in \eqref{4.10} and integrate over $\mathbb{R} \times [0,t]$ to get
		\begin{equation}
			\begin{aligned}
				& \frac{1}{2} \int_{-2L}^{2L} \kappa^{2} \frac{q^{2}}{V^{2}} + \mu \int_{0}^{t} \int_{-2L}^{2L} \kappa^{2} \frac{q_{z}^{2}}{V^{2}(V+r)}
				\\
				= & \frac{1}{2} \int_{-2L}^{2L} \kappa^{2} \frac{q_{0}^{2}}{V^{2}} + s \int_{0}^{t} \int_{-2L}^{2L} \kappa^{2} \frac{V_{z} q^{2}}{V^{3}} + 2 \mu \int_{0}^{t} \int_{-2L}^{2L} \kappa^{2} \frac{V_{z} q q_{z}}{V^{3}(V+r)}
				\\
				& - \int_{0}^{t} \int_{-2L}^{2L} \kappa^{2} \bigg[\frac{s \mu V_zr}{V(V+r)} - (p(V+r)-p(V))\bigg] \left(\frac{q_{z}}{V^{2}} - \frac{2 V_{z} q}{V^{3}}\right)
				\\
				& - \int_{0}^{t} \int_{-2L}^{2L} 2 \kappa \kappa_{z} \left[\mu \left(\frac{q_z}{V+r} + \frac{sV_zr}{V(V+r)}\right)\frac{q}{V^{2}} - (p(V+r)-p(V)) \frac{q}{V^{2}} + \frac{s}{2} \frac{q^{2}}{V^{2}}\right].
			\end{aligned}
			\label{4.7111}
		\end{equation}
		By virtue of \eqref{3-16}, \eqref{4.711}, Taylor's expansion, Young's inequality and the bound $|\kappa_{z}| \leq C/L$, the last term of \eqref{4.7111} can be estimated as
		\begin{align*}
			& \bigg| \int_{0}^{t} \int_{-2L}^{2L} 2 \kappa \kappa_{z} \left[\mu \left(\frac{q_z}{V+r} + \frac{sV_zr}{V(V+r)}\right)\frac{q}{V^{2}} - (p(V+r)-p(V)) \frac{q}{V^{2}} + \frac{s}{2} \frac{q^{2}}{V^{2}}\right] \bigg|
			\\
			\leq & \frac{C}{L} \int_{0}^{t} \int_{-2L}^{2L} \kappa \left( \frac{|q| |q_{z}|}{V^{3}} +  \frac{|V_{z}||r| |q|}{V^{4}} + \frac{|r| |q|}{V^{2}} + \frac{q^{2}}{V^{2}} \right)
			\\
			\leq & \frac{\mu}{6} \int_0^t \int_{-2L}^{2L} \kappa ^{2} \frac{q_z^2}{V^3}+\frac{C}{L^{2}} \int_{0}^{t} \int_{-2L}^{2L} \frac{q^{2}}{V^{3}}+\frac{C}{L}\int_{0}^{t} \int_{-2L}^{2L} \frac{q^{2}+r^{2}}{V^{2}}.
		\end{align*}
		Moreover, Theorem \ref{theorem1.1} gives $V(z) \sim (a_{+} z)^{-1}$ as $z \to +\infty$. Hence, for $z \in [-2L,2L]$, we have $1/V \leq CL$. Using $\|r\|_{w_{2}}^2,\|q\|_{w_{2}}^2 \leq N^2(t)$ and Lemma \ref{lemma3.3}, we obtain
		\begin{equation}
			\begin{aligned}
				& \bigg| \int_{0}^{t} \int_{-2L}^{2L} 2 \kappa \kappa_{z} \left[\mu \left(\frac{q_z}{V+r} + \frac{sV_zr}{V(V+r)}\right)\frac{q}{V^{2}} - (p(V+r)-p(V)) \frac{q}{V^{2}} + \frac{s}{2} \frac{q^{2}}{V^{2}}\right] \bigg|
				\\
				\leq & \frac{\mu}{6} \int_0^t \int_{-2L}^{2L} \kappa^{2} \frac{q_z^2}{V^3} + \frac{C T N^2(T)}{L}.
			\end{aligned}
			\label{3-371}
		\end{equation}
		Substituting \eqref{3-371} into \eqref{4.7111} and noting \eqref{4.711}, we deduce that
		\begin{align*}
			& \frac{1}{2} \int_{-2L}^{2L} \kappa^{2} \frac{q^{2}}{V^{2}} + \left(\frac{2\mu}{3} - \frac{\mu}{6}\right) \int_{0}^{t} \int_{-2L}^{2L} \kappa^{2} \frac{q_{z}^{2}}{V^{3}}
			\\
			\leq & \frac{1}{2} \int_{-2L}^{2L} \kappa^{2} \frac{q_{0}^{2}}{V^{2}} + \frac{C T N^2(T)}{L} + s \int_{0}^{t} \int_{-2L}^{2L} \kappa^{2} \frac{V_{z} q^{2}}{V^{3}} + 2 \mu \int_{0}^{t} \int_{-2L}^{2L} \kappa^{2} \frac{V_{z} q q_{z}}{V^{3}(V+r)}
			\\
			& - \int_{0}^{t} \int_{-2L}^{2L} \kappa^{2} \bigg[\frac{s \mu V_zr}{V(V+r)} - (p(V+r)-p(V))\bigg] \left(\frac{q_{z}}{V^{2}} - \frac{2 V_{z} q}{V^{3}}\right).
		\end{align*}
		Letting $L \to +\infty$, we obtain
		\begin{align*}
			& \frac{1}{2} \int \frac{q^{2}}{V^{2}} + \frac{\mu}{2} \int_{0}^{t} \int \frac{q_{z}^{2}}{V^{3}}
			\\
			\leq & C \int \frac{q_{0}^{2}}{V^{2}} + \int_{0}^{t} \int \left|\frac{s V_{z} q^{2}}{V^{3}}\right| + \int_{0}^{t} \int \left|\frac{2 \mu V_{z} q q_{z}}{V^{3}(V+r)}\right|
			\\
			& + \int_{0}^{t} \int \left|\bigg[\frac{s \mu V_zr}{V(V+r)} - (p(V+r)-p(V))\bigg] \left(\frac{q_{z}}{V^{2}} - \frac{2 V_{z} q}{V^{3}}\right)\right|.
		\end{align*}
		Thanks to $|V_{z}| \leq C V^{2}$ and $2/(3V) \leq 1/(V+r) \leq 2/V$, using Taylor's expansion and Young's inequality, we get
		\begin{equation*}
			\begin{aligned}
				& \int_{0}^{t} \int \left|\frac{s V_{z} q^{2}}{V^{3}}\right| + \int_{0}^{t} \int \left|\frac{2 \mu V_{z} q q_{z}}{V^{3}(V+r)}\right| + \int_{0}^{t} \int \left|\bigg[\frac{s \mu V_zr}{V(V+r)} - (p(V+r)-p(V))\bigg] \left(\frac{q_{z}}{V^{2}} - \frac{2 V_{z} q}{V^{3}}\right)\right|
				\\
				\leq & C \int_{0}^{t} \int \frac{q^{2}}{V} + C \int_{0}^{t} \int\frac{|q| |q_z|}{V^2} + C \int_{0}^{t} \int\frac{|r| |q_z|}{V^2} + C \int_{0}^{t} \int\frac{|r| |q|}{V}
				\\
				\leq & \frac{\mu}{4} \int_{0}^{t} \int\frac{q_z^2}{V^3} + C \int_{0}^{t} \int\frac{q^2+r^2}{V}.
			\end{aligned}
		\end{equation*}
		Combining the above estimates with Lemmas \ref{lemma3.3} and \ref{lemma3.4}, we obtain
		\begin{equation}
			\begin{aligned}
				\frac{1}{2} \int \frac{q^{2}}{V^{2}} + \frac{\mu}{4} \int_{0}^{t} \int \frac{q_{z}^{2}}{V^{3}} \leq & C \int \frac{q_{0}^{2}}{V^{2}} + C \int_{0}^{t} \int \frac{q^{2}}{V} + C \int_{0}^{t} \int \frac{r^{2}}{V}
				\\
				\leq & C (\|\phi_0\|^2 + \|\psi_0\|^2 + \|r_{0}\|_{w_2}^2 +  \|q_{0}\|_{w_{2}}^2).
			\end{aligned}
			\label{3-47}
		\end{equation}
		The proof of this lemma is now complete.
	\end{proof}
	
	\begin{lemma}\label{lemma3.5}
		Let the assumptions of Theorem \ref{theorem3.1} hold. If $N(T) \ll 1$, then there exists a constant $C>0$ independent of $T$ such that
		\begin{equation}
			\|r(t)\|^2+\|\psi_z(t)\|^2
			+\int_0^t \|r_z\|_{w_1}^2 \mathrm{d}\tau + \int_0^t \int d(V)\psi_z^2 \mathrm{d}\tau
			\leq C(\|\phi_0\|_1^2+\|\psi_0\|_1^2),
			\label{3-41}
		\end{equation}
		for any $t \in[0, T]$.
	\end{lemma}
	\begin{proof}
		Differentiating \eqref{3-6} with respect to $z$ gives
		\begin{equation}
			\begin{cases}
				r_t = sr_{z} + \psi_{zz} + \mu(r/V)_{zz} + \mu W_{zz},
				\\
				\psi_{zt} = s\psi_{zz} - p''(V) V_{z}r - p'(V)r_{z} - G_{z}.
			\end{cases}
			\label{3-42}
		\end{equation}
		Multiplying $\eqref{3-42}_{1}$ by $r$ and $\eqref{3-42}_{2}$ by $(-\frac{\psi_z}{p'(V)})$, respectively, and then adding them, we derive
		\begin{equation}
			\begin{aligned}
				\frac{1}{2} \frac{\mathrm{d}}{\mathrm{d}t} \left(r^2 - \frac{\psi_z^2}{p'(V)}\right) = & \left[\mu \left(\frac{r_{z}}{V} - \frac{V_{z} r}{V^{2}}\right) r + \frac{s}{2} r^{2} + \psi_{z} r + \mu W_{z} r + \frac{s\psi_z^2}{-2p'(V)}\right]_{z} - \mu \frac{r_{z}^2}{V}
				\\
				& + \mu \frac{V_z}{V^2} r r_{z} - \frac{s p''(V) V_z}{2 p'(V)^2} \psi_z^2 - \mu W_{z} r_{z} + \frac{p''(V) V_{z}}{p'(V)} r \psi_z +  G_z \frac{\psi_z}{p'(V)}.
			\end{aligned}
			\label{3-4211}
		\end{equation}
		We multiply \eqref{3-4211} by the cut--off function $\kappa^{2}(z)$ and integrate over $\mathbb{R} \times [0,t]$ to get
		\begin{equation}
			\begin{aligned}
				& \frac{1}{2} \int_{-2L}^{2L} \kappa^{2} \left(r^2 - \frac{\psi_z^2}{p'(V)}\right) + \mu \int_{0}^{t} \int_{-2L}^{2L} \kappa^{2} \frac{r_{z}^2}{V} + \int_{0}^{t} \int_{-2L}^{2L} \kappa^{2} \frac{s p''(V) V_z}{2 p'(V)^2} \psi_z^2
				\\
				= & \frac{1}{2} \int_{-2L}^{2L} \kappa^{2} \left(r_{0}^2 - \frac{\psi_{0z}^2}{p'(V)}\right) + \int_{0}^{t} \int_{-2L}^{2L} \kappa^{2} \left(\mu \frac{V_z}{V^2} r r_{z} - \mu W_{z} r_{z} + \frac{p''(V) V_{z}}{p'(V)} r \psi_z + G_z \frac{\psi_z}{p'(V)}\right)
				\\
				& - \int_{0}^{t} \int_{-2L}^{2L} 2 \kappa \kappa_{z} \bigg[\mu \left(\frac{r_{z}}{V} - \frac{V_{z} r}{V^{2}}\right) r + \frac{s}{2} r^{2} + \psi_{z} r + \mu W_{z} r + \frac{s\psi_z^2}{-2p'(V)}\bigg].
			\end{aligned}
			\label{3-4311}
		\end{equation}
		For the last term of \eqref{3-4311}, the estimates $|\kappa_{z}| \leq C/L$, $\|\psi_{z}\|, \|r/V\|_{L^{\infty}}, \|r_{z}/V\|_{L^{\infty}} \leq C N(t)$, together with \eqref{3-18} and Lemma \ref{lemma3.3}, imply that
		\begin{align*}
			& \left|\int_{0}^{t} \int_{-2L}^{2L} 2 \kappa \kappa_{z} \bigg[\mu \left(\frac{r_{z}}{V} - \frac{V_{z} r}{V^{2}}\right) r + \frac{s}{2} r^{2} + \psi_{z} r + \mu W_{z} r + \frac{s\psi_z^2}{-2p'(V)}\bigg]\right|
			\\
			\leq & \frac{C}{L} \int_{0}^{t} \int_{-2L}^{2L} \kappa \bigg(\frac{|r_{z}| |r|}{V} + r^{2} + |\psi_{z}| |r| + \frac{|r|^{3}}{V} + \frac{|r_{z}| |r|^{2}}{V^{2}} + \psi_{z}^{2} \bigg)
			\\
			\leq & \frac{C (1+N(T))}{L} \int_{0}^{t} \int_{-2L}^{2L} \kappa^{2} \frac{r_{z}^{2}}{V} + \frac{C(1+N(T))}{L} \int_{0}^{t} \int_{-2L}^{2L} \frac{r^{2}}{V} + \frac{C}{L} \int_{0}^{t} \int_{-2L}^{2L} \psi_{z}^{2}
			\\
			\leq & \frac{C (1+N(T))}{L} \int_{0}^{t} \int_{-2L}^{2L} \kappa^{2} \frac{r_{z}^{2}}{V} + \frac{C (1+N(T))}{L} (\|\phi_0\|^2 + \|\psi_0\|^2) + \frac{C T N^2(T)}{L}.
		\end{align*}
		Therefore, \eqref{3-4311} becomes
		\begin{align*}
			& \frac{1}{2} \int_{-2L}^{2L} \kappa^{2} \left(r^2 - \frac{\psi_z^2}{p'(V)}\right) + \left(\mu - \frac{C (1+N(T))}{L}\right) \int_{0}^{t} \int_{-2L}^{2L} \kappa^{2} \frac{r_{z}^2}{V} + \int_{0}^{t} \int_{-2L}^{2L} \kappa^{2} \frac{s p''(V) V_z}{2 p'(V)^2} \psi_z^2
			\\
			\leq & \frac{1}{2} \int_{-2L}^{2L} \kappa^{2} \left(r_{0}^2 - \frac{\psi_{0z}^2}{p'(V)}\right) + \frac{C (1+N(T))}{L} (\|\phi_0\|^2 + \|\psi_0\|^2) + \frac{C T N^{2}(T)}{L}
			\\
			& + \int_{0}^{t} \int_{-2L}^{2L} \kappa^{2} \left(\mu \frac{V_z}{V^2} r r_{z} - \mu W_{z} r_{z} + \frac{p''(V) V_{z}}{p'(V)} r \psi_z + G_z \frac{\psi_z}{p'(V)}\right).
		\end{align*}
		Passing to the limit $L \to + \infty$, we obtain
		\begin{equation}
			\begin{aligned}
				& \frac{1}{2} \int \left(r^2 - \frac{\psi_z^2}{p'(V)}\right) + \mu \int_{0}^{t} \int \frac{r_{z}^2}{V} + \int_{0}^{t} \int \frac{s p''(V) V_z}{2 p'(V)^2} \psi_z^2
				\\
				\leq & C (\|r_{0}\|^2 + \|\psi_{0z}\|^2) + \int_{0}^{t} \int \left|\mu \frac{V_z}{V^2} r r_{z} - \mu W_{z} r_{z} + \frac{p''(V) V_{z}}{p'(V)} r \psi_z + G_z \frac{\psi_z}{p'(V)}\right|.
			\end{aligned}
			\label{3-4411}
		\end{equation}
		Employing \eqref{3-16}, \eqref{3-18}, $\|r_{z}/V\|_{L^\infty} \leq CN(t)$, $\|r/V\|_{L^\infty}\leq CN(t)$, together with Young's inequality, we have
		\begin{equation}
			\begin{aligned}
				& \int_{0}^{t} \int \left|\mu \frac{V_z}{V^2} r r_{z} - \mu W_{z} r_{z} + \frac{p''(V) V_{z}}{p'(V)} r \psi_z + G_z \frac{\psi_z}{p'(V)}\right|
				\\
				\leq & C \int_{0}^{t} \int \left(|r r_{z}| + |V_{z}| |r \psi_z| + \frac{r^2 |r_{z}|}{V} + \frac{|r| r_{z}^2}{V^2} + |V_z| r^2 |\psi_z| + |r| |r_{z}| |\psi_z|\right)
				\\
				\leq & \left(\frac{\mu}{8} + C N(t)\right) \int_{0}^{t} \int \frac{r_{z}^2}{V} + C \int_{0}^{t} \int \frac{r^2}{V} + \frac{1}{4} \int_{0}^{t} \int \frac{s p''(V) V_z}{2p'(V)^2} \psi_z^2.
			\end{aligned}
			\label{3-44}
		\end{equation}
		Combining \eqref{3-22} and \eqref{3-44} with \eqref{3-4411}, we conclude that
		\begin{align*}
			\int \left(r^2 - \frac{\psi_z^2}{p'(V)}\right) + \int_0^t \int \frac{r_{z}^2}{V} + \int_0^t \int d(V) \psi_z^2 \leq C \left(\int r_{0}^2+\int \psi_{0 z}^2+\int \phi_0^2+\int \psi_0^2\right).
		\end{align*}
		We thus finish the proof of \eqref{3-41}.
	\end{proof}
	
	\begin{lemma}\label{lemma3.8}
		Let the assumptions of Theorem \ref{theorem3.1} hold. If $N(T) \ll 1$, then there exists a constant $C>0$ independent of $T$ such that
		\begin{equation}
			\|q_z(t)\|_{w_2}^2+\int_0^t\|q_{zz}\|_{w_3}^2\,\mathrm{d}\tau
			\leq C (\|\phi_0\|^2 + \|\psi_0\|_{1}^2 + \|r_{0}\|_{w_2}^2 +  \|q_{0}\|_{1,w_{2}}^2),
			\label{3-53}
		\end{equation}
		for any $t \in[0, T]$.
	\end{lemma}
	\begin{proof}
		Differentiating \eqref{4.7} gives
		\begin{equation}
			\begin{aligned}
				q_{zt} = & sq_{zz} + \left[\mu\left(\frac{q_z}{V+r} + \frac{sV_zr}{V(V+r)}\right) - (p(V+r)-p(V))\right]_{zz}
				\\
				= & sq_{zz} + \bigg\{\mu\left[\frac{q_{zz}}{V+r} - \frac{(V_{z} + r_{z})q_z}{(V+r)^{2}} + s\left(\frac{V_{z}}{V^{2}}\right)_{z} \frac{V r}{V+r} + s\frac{V_{z}}{V^{2}} \left(\frac{V r}{V+r}\right)_{z}\right]
				\\
				& \quad - (p(V+r)-p(V))_{z}\bigg\}_{z}.
			\end{aligned}
			\label{4.8}
		\end{equation}
		Multiplying \eqref{4.8} by $\frac{q_{z}}{V^{2}}$, we get
		\begin{equation}
			\begin{aligned}
				\frac{1}{2} \frac{\mathrm{d}}{\mathrm{d}t} \frac{q_{z}^{2}}{V^{2}} = & \left\{\mu \left[\frac{q_{zz}}{V+r} - \frac{(V_{z} + r_{z})q_z}{(V+r)^{2}} + s\left(\frac{V_{z}}{V^{2}}\right)_{z} \frac{V r}{V+r} + s\frac{V_{z}}{V^{2}} \left(\frac{V r}{V+r}\right)_{z}\right] \frac{q_{z}}{V^{2}} \right.
				\\
				& \left. - (p(V + r) - p(V))_{z} \frac{q_{z}}{V^{2}} + \frac{s}{2} \frac{q_{z}^{2}}{V^{2}} \right\}_{z} + s \frac{V_{z} q_{z}^{2}}{V^{3}}
				\\
				& - \mu \left[\frac{q_{zz}}{V+r} - \frac{(V_{z} + r_{z})q_z}{(V+r)^{2}} + s\left(\frac{V_{z}}{V^{2}}\right)_{z} \frac{V r}{V+r} + s\frac{V_{z}}{V^{2}} \left(\frac{V r}{V+r}\right)_{z}\right]
				\\
				& \, \times \left(\frac{q_{zz}}{V^2}-\frac{2V_{z}q_z}{V^3}\right) + [p(V + r) - p(V)]_{z} \left(\frac{q_{zz}}{V^2}-\frac{2V_{z}q_z}{V^3}\right).
			\end{aligned}
			\label{4.9}
		\end{equation}
		Then, multiplying \eqref{4.9} by $\kappa^{2}(z)$  and integrating over $\mathbb{R} \times [0,t]$, we obtain
		\begin{equation}
			\begin{aligned}
				& \frac{1}{2} \int_{-2L}^{2L} \kappa^{2} \frac{q_{z}^{2}}{V^{2}} + \mu \int_{0}^{t} \int_{-2L}^{2L} \frac{\kappa^{2} q_{zz}^{2}}{V^{2}(V+r)}
				\\
				= & \frac{1}{2} \int_{-2L}^{2L} \kappa^{2} \frac{q_{0z}^{2}}{V^{2}} + s \int_{0}^{t} \int_{-2L}^{2L} \kappa^{2} \frac{V_{z} q_{z}^{2}}{V^{3}} + 2\mu \int_{0}^{t} \int_{-2L}^{2L} \frac{\kappa^{2} V_{z} q_{z z} q_{z}}{V^{3}(V+r)}
				\\
				& + \mu \int_{0}^{t} \int_{-2L}^{2L} \kappa^{2} \frac{q_{z z} q_{z} (V_{z} + r_{z})}{V^{2} (V+r)^{2}} - 2 \mu \int_{0}^{t} \int_{-2L}^{2L} \kappa^{2} \frac{q_{z}^{2} V_{z} (V_{z} + r_{z})}{V^{3} (V+r)^{2}}
				\\
				& + \int_{0}^{t} \int_{-2L}^{2L} \kappa^{2} \left[(p(V + r) - p(V))_{z} - s \mu \left(\frac{V_{z}}{V^{2}}\right)_{z} \frac{V r}{V+r} - s \mu \frac{V_{z}}{V^{2}} \left(\frac{V r}{V+r}\right)_{z}\right]
				\\
				& \times \left(\frac{q_{zz}}{V^2}-\frac{2V_{z}q_z}{V^3}\right) - \int_{0}^{t} \int_{-2L}^{2L} 2 \kappa \kappa_{z} \left\{\mu \left[\frac{q_{zz}}{V+r} - \frac{(V_{z} + r_{z})q_z}{(V+r)^{2}} + s\left(\frac{V_{z}}{V^{2}}\right)_{z} \frac{V r}{V+r} \right. \right.
				\\
				& \left. \left.+ s\frac{V_{z}}{V^{2}} \left(\frac{V r}{V+r}\right)_{z}\right] \frac{q_{z}}{V^{2}} - (p(V + r) - p(V))_{z} \frac{q_{z}}{V^{2}} + \frac{s}{2} \frac{q_{z}^{2}}{V^{2}} \right\}.
			\end{aligned}
			\label{4.11}
		\end{equation}
		Now, we estimate the last term of \eqref{4.11}. Applying Young's inequality, Lemmas \ref{lemma3.3}, \ref{lemma3.6} and \ref{lemma3.5}, together with the estimates \eqref{3-16}, $|(V_{z}/V^{2})_{z}| \leq C |V_{z}|$, $\|r_{z}/V\|_{L^{\infty}} \leq C N(T)$, $2/(3V) \leq 1/(V+r) \leq 2/V$ and \eqref{4.10}, we have
		\begin{align*}
			& \bigg| \int_{0}^{t} \int_{-2L}^{2L} 2 \kappa \kappa_{z} \left\{\mu \left[\frac{q_{zz}}{V+r} - \frac{(V_{z} + r_{z})q_z}{(V+r)^{2}} + s\left(\frac{V_{z}}{V^{2}}\right)_{z} \frac{V r}{V+r} \right. \right.
			\\
			& \left. \left.+ s\frac{V_{z}}{V^{2}} \left(\frac{V r}{V+r}\right)_{z}\right] \frac{q_{z}}{V^{2}} - (p(V + r) - p(V))_{z} \frac{q_{z}}{V^{2}} + \frac{s}{2} \frac{q_{z}^{2}}{V^{2}} \right\} \bigg|
			\\
			\leq & \frac{C}{L} \int_{0}^{t} \int_{-2L}^{2L} \kappa \left(\frac{|q_{z}||q_{zz}|}{V^{3}} + \frac{|V_{z}|q_{z}^{2}}{V^{4}} + \frac{|r_{z}| q_{z}^{2}}{V^{4}} + \frac{|V_{z}| |r| |q_{z}|}{V^{2}} + \frac{|r_{z}| |q_{z}|}{V^{2}} \right.
			\\
			& \left. \quad + \frac{|V_{z}| |r| |q_{z}|}{V^{3}} + \frac{|r_{z}| |r| |q_{z}|}{V^{3}} + \frac{q_{z}^{2}}{V^{2}}\right)
			\\
			\leq & \frac{C}{L} \int_{0}^{t} \int_{-2L}^{2L} \left( \kappa^{2} \frac{q_{zz}^{2}}{V^{3}} + \frac{q_{z}^{2}}{V^{3}} + \frac{r^{2}}{V} + \frac{r_{z}^{2}}{V} \right)
			\\
			\leq & \frac{C}{L} \int_{0}^{t} \int_{-2L}^{2L} \kappa^{2} \frac{q_{zz}^{2}}{V^{3}} + \frac{C}{L} (\|\phi_0\|^2 + \|\psi_{0}\|^2 + \|r_{0}\|_{w_2}^2 + \|q_{0}\|_{w_{2}}^2).
		\end{align*}
		Based on the above inequality and $2/(3V) \leq (V+r)^{-1}$, \eqref{4.11} becomes
		\begin{align*}
			& \frac{1}{2} \int_{-2L}^{2L} \kappa^{2} \frac{q_{z}^{2}}{V^{2}} + \left(\frac{2\mu}{3} - \frac{C}{L}\right)  \int_{0}^{t} \int_{-2L}^{2L} \kappa^{2} \frac{q_{zz}^{2}}{V^{3}}
			\\
			\leq & \frac{1}{2} \int_{-2L}^{2L} \kappa^{2} \frac{q_{0z}^{2}}{V^{2}} + \frac{C}{L} (\|\phi_0\|^2 + \|\psi_{0}\|^2 + \|r_{0}\|_{w_2}^2 + \|q_{0}\|_{w_{2}}^2)+ s \int_{0}^{t} \int_{-2L}^{2L} \kappa^{2} \frac{V_{z} q_{z}^{2}}{V^{3}}
			\\
			& + 2\mu \int_{0}^{t} \int_{-2L}^{2L} \frac{\kappa^{2} V_{z} q_{z z} q_{z}}{V^{3}(V+r)} + \mu \int_{0}^{t} \int_{-2L}^{2L} \kappa^{2} \frac{q_{z z} q_{z} (V_{z} + r_{z})}{V^{2} (V+r)^{2}} - 2 \mu \int_{0}^{t} \int_{-2L}^{2L} \kappa^{2} \frac{q_{z}^{2} V_{z} (V_{z} + r_{z})}{V^{3} (V+r)^{2}}
			\\
			& + \int_{0}^{t} \int_{-2L}^{2L} \kappa^{2} \left[(p(V + r) - p(V))_{z} - s \mu \left(\frac{V_{z}}{V^{2}}\right)_{z} \frac{V r}{V+r} - s \mu \frac{V_{z}}{V^{2}} \left(\frac{V r}{V+r}\right)_{z}\right] \left(\frac{q_{zz}}{V^2}-\frac{2V_{z}q_z}{V^3}\right).
		\end{align*}
		Letting $L \to +\infty$ and then using $(V+r)^{-1} \leq 2/V$, we have
		\begin{equation}
			\begin{aligned}
				& \frac{1}{2} \int \frac{q_{z}^{2}}{V^{2}} + \frac{\mu}{2} \int_{0}^{t} \int \frac{q_{zz}^{2}}{V^{3}}
				\\
				\leq & C \int \frac{q_{0z}^{2}}{V^{2}} + C \int_{0}^{t} \int \frac{|V_{z}| q_{z}^{2}}{V^{3}} + C \int_{0}^{t} \int \frac{|V_{z}| |q_{z z}| |q_{z}|}{V^{4}} + C \int_{0}^{t} \int \frac{|r_{z}| |q_{z z}| |q_{z}|}{V^{4}}
				\\
				& + C \int_{0}^{t} \int \frac{q_{z}^{2} V_{z}^{2}}{V^{5}} + C \int_{0}^{t} \int \frac{|V_{z}| q_{z}^{2} |r_{z}|}{V^{5}} + C \int_{0}^{t} \int \bigg|(p(V + r) - p(V))_{z}
				\\
				& \quad - s \mu \left(\frac{V_{z}}{V^{2}}\right)_{z} \frac{V r}{V+r} - s \mu \frac{V_{z}}{V^{2}} \left(\frac{V r}{V+r}\right)_{z}\bigg| \left(\frac{|q_{zz}|}{V^2}+\frac{2|V_{z}||q_z|}{V^3}\right).
			\end{aligned}
			\label{4.13}
		\end{equation}
		Combining the bound $1 \leq C/V$ with \eqref{3-16}, \eqref{profile:coeff} and \eqref{rz:infinity} and applying Lemmas \ref{lemma3.3}, \ref{lemma3.6} and \ref{lemma3.5}, we obtain
		\begin{align*}
			& \int_{0}^{t} \int \frac{|V_{z}| q_{z}^{2}}{V^{3}} + \int_{0}^{t} \int \frac{|V_{z}| |q_{z z}| |q_{z}|}{V^{4}} + \int_{0}^{t} \int \frac{|r_{z}| |q_{z z}| |q_{z}|}{V^{4}} + \int_{0}^{t} \int \frac{q_{z}^{2} V_{z}^{2}}{V^{5}} + \int_{0}^{t} \int \frac{|V_{z}| q_{z}^{2} |r_{z}|}{V^{5}}
			\\
			\leq & C \int_{0}^{t} \int \frac{q_{z}^{2}}{V} + \left(\frac{\mu}{8} + CN(t)\right) \int_{0}^{t} \int \frac{q_{zz}^{2}}{V^{3}} + C N(t) \int_{0}^{t} \int \left(\frac{q_{z}^{2}}{V^{3}} + \frac{q_{z}^{2}}{V^{2}}\right)
			\\
			\leq & C (\|\phi_0\|^2 + \|\psi_0\|_{1}^2 + \|r_{0}\|_{w_2}^2 + \|q_{0}\|_{w_{2}}^2) + \left(\frac{\mu}{8} + CN(t)\right) \int_{0}^{t} \int \frac{q_{zz}^{2}}{V^{3}},
		\end{align*}
		and
		\begin{align*}
			& \int_{0}^{t} \int \bigg|(p(V + r) - p(V))_{z} - s \mu \left(\frac{V_{z}}{V^{2}}\right)_{z} \frac{V r}{V+r} - s \mu \frac{V_{z}}{V^{2}} \left(\frac{V r}{V+r}\right)_{z}\bigg| \left(\frac{|q_{zz}|}{V^2}+\frac{2|V_{z}||q_z|}{V^3}\right)
			\\
			\leq & C \int_{0}^{t} \int \bigg(|r_{z}| + |V_{z}| |r| + V |r| + \frac{|r| |r_{z}|}{V}\bigg) \left(\frac{|q_{zz}|}{V^2}+\frac{2|V_{z}||q_z|}{V^3}\right)
			\\
			\leq & C \int_{0}^{t} \int \frac{r^{2}}{V} + C (1+N(t)) \int_{0}^{t} \int \frac{r_{z}^{2}}{V} + \left(\frac{\mu}{8} + CN(t)\right) \int_{0}^{t} \int \frac{q_{zz}^{2}}{V^{3}} + C (1+N(t)) \int_{0}^{t} \int \frac{q_{z}^{2}}{V}
			\\
			\leq & C (\|\phi_0\|^2 + \|\psi_0\|_{1}^2 + \|r_{0}\|_{w_2}^2 + \|q_{0}\|_{w_{2}}^2) + \left(\frac{\mu}{8} + CN(t)\right) \int_{0}^{t} \int \frac{q_{zz}^{2}}{V^{3}}.
		\end{align*}
		In summary, \eqref{4.13} takes the form
		\begin{equation*}
			\frac{1}{2} \int \frac{q_{z}^{2}}{V^{2}} + \frac{\mu}{8} \int_{0}^{t} \int \frac{q_{zz}^{2}}{V^{3}} \leq C (\|\phi_0\|^2 + \|\psi_0\|_{1}^2 + \|r_{0}\|_{w_2}^2 + \|q_{0}\|_{1, w_{2}}^2),
		\end{equation*}
		provided $N(t)$ is suitably small.
	\end{proof}
	
	\begin{lemma}\label{lemma3.9}
		Let the assumptions of Theorem \ref{theorem3.1} hold. If $N(T) \ll 1$, then there exists a constant $C>0$ independent of $T$ such that
		\begin{equation}
			\|r_z(t)\|_{w_2}^2 \leq C (\|\phi_0\|^2 + \|\psi_0\|_{1}^2 + \|r_{0}\|_{1,w_2}^2 + \|q_{0}\|_{1, w_{2}}^2),
			\label{3-57}
		\end{equation}
		for any $t \in[0, T]$.
	\end{lemma}
	\begin{proof}
		Differentiating $r_t=s r_z+q_z$ with respect to $z$, multiplying the resulting equation by $\kappa^{2} r_z/V^2$ and integrating over $\mathbb{R} \times [0,t]$, we obtain
		\begin{equation*}
			\frac{1}{2} \int_{-2L}^{2L} \kappa^{2} \frac{r_z^2}{V^2} = \frac{1}{2} \int_{-2L}^{2L} \kappa^{2} \frac{r_{0z}^{2}}{V^{2}} + \int_{0}^{t} \int_{-2L}^{2L} \kappa^{2} \frac{q_{zz}r_z}{V^2} + s \int_{0}^{t} \int_{-2L}^{2L} \kappa^{2} \frac{V_zr_z^2}{V^3} - s \int_{0}^{t} \int_{-2L}^{2L} \kappa \kappa_{z}  \frac{r_{z}^{2}}{V^{2}}.
		\end{equation*}
		By \eqref{3-11}, we have
		\begin{equation*}
			\left| - s \int_{0}^{t} \int_{-2L}^{2L} \kappa \kappa_{z}  \frac{r_{z}^{2}}{V^{2}} \right| \leq \frac{CTN^{2}(T)}{L}.
		\end{equation*}
		Letting $L \to + \infty$ and using Young's inequality, we get
		\begin{equation*}
			\frac{1}{2} \int \frac{r_z^2}{V^2} \leq \frac{1}{2} \int \frac{r_{0z}^{2}}{V^{2}} + \int_{0}^{t} \int \frac{q_{zz}r_z}{V^2} + s \int_{0}^{t} \int \frac{V_zr_z^2}{V^3} \leq \frac{1}{2} \int \frac{r_{0z}^{2}}{V^{2}} + C \int_{0}^{t} \int \frac{q_{zz}^{2}}{V^{3}} + C \int_{0}^{t} \int \frac{r_{z}^{2}}{V}.
		\end{equation*}
		Applying Lemmas \ref{lemma3.5} and \ref{lemma3.8}, we complete the proof.
	\end{proof}
	
	\begin{lemma}\label{lemma3.11}
		Let the assumptions of Theorem \ref{theorem3.1} hold. If $N(T) \ll 1$, then there exists a constant $C>0$ independent of $T$ such that
		\begin{equation}
			\|r_{z}\|^2 + \|\psi_{zz}\|^2 + \int_0^t \|r_{zz}\|_{w_1}^2  + \int_0^t \int d(V) \psi_{zz}^2 \leq C (\|\phi_0\|_2^2 + \|\psi_0\|_2^2),
			\label{3-91}
		\end{equation}
		for any $t \in [0, T]$.
	\end{lemma}
	\begin{proof}
		Differentiating \eqref{3-42} with respect to $z$ gives
		\begin{equation}
			\begin{cases}
				r_{zt} = sr_{zz} + \psi_{zzz} + \mu(r/V)_{zzz} + \mu W_{zzz},
				\\
				\psi_{zzt} = s\psi_{zzz} - p'(V) r_{zz} - 2p''(V)V_zr_z - p'''(V)V_z^2 r - p''(V)V_{zz}r - G_{zz}.
			\end{cases}
			\label{3-421}
		\end{equation}
		Multiplying $\eqref{3-421}_{1}$ by $r_{z}$ and $\eqref{3-421}_{2}$ by $(-\frac{\psi_{zz}}{p'(V)})$, respectively, and then adding them yields
		\begin{equation}
			\begin{aligned}
				& \frac{1}{2} \frac{\mathrm{d}}{\mathrm{d}t} \left(r_{z}^2 - \frac{\psi_{zz}^2}{p'(V)}\right)
				\\
				= & \left[\mu \left(\frac{r_{zz}}{V} - \frac{2V_{z} r_{z}}{V^{2}} - \frac{V_{zz} r}{V^{2}} + \frac{2V_{z}^{2} r}{V^{3}}\right) r_{z} + \frac{s}{2} r_{z}^{2} + \psi_{zz} r_{z} + \mu W_{zz} r_{z} + \frac{s\psi_{zz}^2}{-2p'(V)}\right]_{z} - \mu \frac{r_{zz}^2}{V}
				\\
				& + 2\mu \frac{V_z}{V^2} r_{z} r_{zz} + \mu \frac{V_{zz}}{V^2} r r_{zz} - 2\mu \frac{V_z^{2}}{V^3} r r_{zz}- \frac{s p''(V) V_z}{2 p'(V)^2} \psi_{zz}^2 - \mu W_{zz} r_{zz}
				\\
				& + \frac{2p''(V)V_z r_z \psi_{zz}}{p'(V)} + \frac{ p'''(V)V_z^2 r \psi_{zz} }{p'(V)} + \frac{ p''(V)V_{zz} r \psi_{zz}}{p'(V)} +  G_{zz} \frac{\psi_{zz}}{p'(V)}.
			\end{aligned}
			\label{3-411}
		\end{equation}
		Multiplying \eqref{3-411} by $\kappa^{2}(z)$ and integrating over $\mathbb{R} \times [0,t]$, we get
		\begin{equation}
			\begin{aligned}
				& \frac{1}{2} \int_{-2L}^{2L} \kappa^{2} \left(r_{z}^2 - \frac{\psi_{zz}^2}{p'(V)}\right) + \mu \int_{0}^{t} \int_{-2L}^{2L} \kappa^{2} \frac{r_{zz}^2}{V} + \int_{0}^{t} \int_{-2L}^{2L} \kappa^{2} \frac{s p''(V) V_z}{2 p'(V)^2} \psi_{zz}^2
				\\
				= & \frac{1}{2} \int_{-2L}^{2L} \kappa^{2} \left(r_{0z}^2 - \frac{\psi_{0zz}^2}{p'(V)}\right) + \int_{0}^{t} \int_{-2L}^{2L} \kappa^{2} \left(2\mu \frac{V_z}{V^2} r_{z} r_{zz} + \mu \frac{V_{zz}}{V^2} r r_{zz} - 2\mu \frac{V_z^{2}}{V^3} r r_{zz} - \mu W_{zz} r_{zz} \right.
				\\
				& \left. \quad + \frac{2p''(V)V_z r_z \psi_{zz}}{p'(V)} + \frac{ p'''(V)V_z^2 r \psi_{zz} }{p'(V)} + \frac{ p''(V)V_{zz} r \psi_{zz}}{p'(V)} +  G_{zz} \frac{\psi_{zz}}{p'(V)}\right) - \int_{0}^{t} \int_{-2L}^{2L} 2 \kappa \kappa_{z}
				\\
				& \quad \times \bigg[\mu \left(\frac{r_{zz}}{V} - \frac{2V_{z} r_{z}}{V^{2}} - \frac{V_{zz} r}{V^{2}} + \frac{2V_{z}^{2} r}{V^{3}}\right) r_{z} + \frac{s}{2} r_{z}^{2} + \psi_{zz} r_{z} + \mu W_{zz} r_{z} + \frac{s\psi_{zz}^2}{-2p'(V)}\bigg].
			\end{aligned}
			\label{3-451}
		\end{equation}
		For the last term of \eqref{3-451}, the estimates $\|\psi_{zz}\| \leq N(t)$, $\|r/V\|_{L^{\infty}},\|r_{z}/V\|_{L^{\infty}} \leq C N(t)$, together with \eqref{3-16}, \eqref{3-18}, Lemmas \ref{lemma3.3} and \ref{lemma3.5}, imply that
		\begin{align*}
			& \left|\int_{0}^{t} \int_{-2L}^{2L} 2 \kappa \kappa_{z} \bigg[\mu \left(\frac{r_{zz}}{V} - \frac{2V_{z} r_{z}}{V^{2}} - \frac{V_{zz} r}{V^{2}} + \frac{2V_{z}^{2} r}{V^{3}}\right) r_{z} + \frac{s}{2} r_{z}^{2} + \psi_{zz} r_{z} + \mu W_{zz} r_{z} + \frac{s\psi_{zz}^2}{-2p'(V)}\bigg]\right|
			\\
			\leq & \frac{C}{L} \int_{0}^{t} \int_{-2L}^{2L} \kappa \bigg(\frac{|r_{zz}| |r_{z}|}{V} + r_{z}^{2} + V|r||r_{z}| + |\psi_{zz}| |r_{z}| + (|V_{zz}|+|r_{z}|) |r|^{2}
			\\
			& \quad + \frac{|r_{z}|^{3}}{V^{2}} + \frac{|r_{z}|^{2} |r|}{V} + \frac{|r||r_{z}||r_{zz}|}{V^{2}} + \psi_{zz}^{2}\bigg)
			\\
			\leq & \frac{C (1+N(T))}{L} \int_{0}^{t} \int_{-2L}^{2L} \kappa^{2} \frac{r_{zz}^{2}}{V} + \frac{C (1+N(T))}{L} \int_{0}^{t} \int_{-2L}^{2L} \left(\frac{r^{2}}{V} + \frac{r_{z}^{2}}{V}\right) + \frac{C}{L} \int_{0}^{t} \int_{-2L}^{2L} \psi_{zz}^{2}
			\\
			\leq & \frac{C (1+N(T))}{L} \int_{0}^{t} \int_{-2L}^{2L} \kappa^{2} \frac{r_{zz}^{2}}{V} + \frac{C (1+N(T))}{L} (\|\phi_0\|_{1}^2 + \|\psi_0\|_{1}^2) + \frac{C T N^{2}(T)}{L}.
		\end{align*}
		Therefore, \eqref{3-451} becomes
		\begin{align*}
			& \frac{1}{2} \int_{-2L}^{2L} \kappa^{2} \left(r_{z}^2 - \frac{\psi_{zz}^2}{p'(V)}\right) + \left(\mu - \frac{C (1+N(T))}{L}\right) \int_{0}^{t} \int_{-2L}^{2L} \kappa^{2} \frac{r_{zz}^2}{V} + \int_{0}^{t} \int_{-2L}^{2L} \kappa^{2} \frac{s p''(V) V_z}{2 p'(V)^2} \psi_{zz}^2
			\\
			\leq & \frac{1}{2} \int_{-2L}^{2L} \kappa^{2} \left(r_{0z}^2 - \frac{\psi_{0zz}^2}{p'(V)}\right) + \frac{C (1+N(T))}{L} (\|\phi_0\|_{1}^2 + \|\psi_0\|_{1}^2) + \frac{C T N^{2}(T)}{L}
			\\
			& + \int_{0}^{t} \int_{-2L}^{2L} \kappa^{2} \left[2\mu \frac{V_z}{V^2} r_{z} r_{zz} + \mu \frac{V_{zz}}{V^2} r r_{zz} - 2\mu \frac{V_z^{2}}{V^3} r r_{zz} - \mu W_{zz} r_{zz} \right.
			\\
			& \left. \quad + \frac{2p''(V)V_z r_z \psi_{zz}}{p'(V)} + \frac{ p'''(V)V_z^2 r \psi_{zz} }{p'(V)} + \frac{ p''(V)V_{zz} r \psi_{zz}}{p'(V)} +  G_{zz} \frac{\psi_{zz}}{p'(V)}\right].
		\end{align*}
		Letting $L \to + \infty$ and using \eqref{3-16}, \eqref{3-15}, \eqref{rz:infinity} and \eqref{3-18}, we obtain
		\begin{equation}
			\begin{aligned}
				& \frac{1}{2} \int \left(r_{z}^2 - \frac{\psi_{zz}^2}{p'(V)}\right) + \mu \int_{0}^{t} \int \frac{r_{zz}^2}{V} + \int_{0}^{t} \int \frac{s p''(V) V_z}{2 p'(V)^2} \psi_{zz}^2
				\\
				\leq & C (\|r_{0z}\|^2 + \|\psi_{0zz}\|^2) + \int_{0}^{t} \int \left|2\mu \frac{V_z}{V^2} r_{z} r_{zz} + \mu \frac{V_{zz}}{V^2} r r_{zz} - 2\mu \frac{V_z^{2}}{V^3} r r_{zz} - \mu W_{zz} r_{zz} \right.
				\\
				& \left. \quad + \frac{2p''(V)V_z r_z \psi_{zz}}{p'(V)} + \frac{ p'''(V)V_z^2 r \psi_{zz} }{p'(V)} + \frac{ p''(V)V_{zz} r \psi_{zz}}{p'(V)}\right| +  \int_{0}^{t} \int \left|G_{zz} \frac{\psi_{zz}}{p'(V)}\right|.
			\end{aligned}
			\label{3.461}
		\end{equation}
		We now estimate the terms on the right--hand of \eqref{3.461}. Using $|V_{z z}| \leq CV|V_z|$, $\|r/V\|_{L^{\infty}} \leq C N(t)$, $\|r_{z}/V\|_{L^{\infty}} \leq C N(t)$, $\|\psi_{zz}\| \leq C N(t)$ and Young's inequality, we have
		\begin{align*}
			& \int_{0}^{t} \int \left|2\mu \frac{V_z}{V^2} r_{z} r_{zz} + \mu \frac{V_{zz}}{V^2} r r_{zz} - 2\mu \frac{V_z^{2}}{V^3} r r_{zz} - \mu W_{zz} r_{zz} \right.
			\\
			& \left. \quad + \frac{2p''(V)V_z r_z \psi_{zz}}{p'(V)} + \frac{ p'''(V)V_z^2 r \psi_{zz} }{p'(V)} + \frac{ p''(V)V_{zz} r \psi_{zz}}{p'(V)} \right|
			\\
			\leq & C \int_{0}^{t} \int \left(|r_{z} r_{zz}| + V|r r_{zz}| + r^{2}| r_{zz}| + \frac{|r_{zz}| r_{z}^2}{V^2} + \frac{|r| |r_{z}||r_{zz}| }{V}+ \frac{|r| r_{zz}^{2}}{V^{2}} + |V_{z}| |r_{z} \psi_{zz}| + V |V_{z}| |r \psi_{zz}|\right)
			\\
			\leq & \left(\frac{\mu}{8} + C N(t)\right) \int_{0}^{t} \int \frac{r_{zz}^2}{V} + C (1+N(t)) \int_{0}^{t} \int \left(\frac{r^2}{V} + \frac{r_{z}^2}{V}\right) + \frac{1}{4} \int_{0}^{t} \int \frac{s p''(V) V_z}{2p'(V)^2} \psi_{zz}^2,
		\end{align*}
		and
		\begin{align*}
			& \int_{0}^{t} \int \left| G_{zz} \frac{\psi_{zz}}{p'(V)} \right|
			\\
			\leq & C \int_{0}^{t} \int \left(V_z^2r^2 |\psi_{zz}| + |V_z||r||r_{z}||\psi_{zz}| + |V_{zz}|r^2 |\psi_{zz}| + r_{z}^2 |\psi_{zz}| + |r||r_{zz}| |\psi_{zz}|\right)
			\\
			\leq & C \int_{0}^{t} \left( \|r\|_{L^{\infty}} \|r\| + \|r\|_{L^{\infty}} \|r_{z}\| + \|r_{z}\|_{L^{\infty}} \|r_{z}\| + \|r\|_{L^{\infty}} \|r_{zz}\|\right) \|\psi_{zz}\|
			\\
			\leq & C N(t) \int_{0}^{t} \left( \|r_{z}\|^{\frac{1}{2}} \|r\|^{\frac{3}{2}} + \|r_{z}\|^{\frac{3}{2}} \|r\|^{\frac{1}{2}} + \|r_{zz}\|^{\frac{1}{2}} \|r_{z}\|^{\frac{3}{2}} + \|r\|^{\frac{1}{2}} \|r_{z}\|^{\frac{1}{2}} \|r_{zz}\|\right)
			\\
			\leq & C N(t) \int_{0}^{t} \left( \|r\|^{2} + \|r_{z}\|^{2} + \|r_{zz}\|^{2}\right).
		\end{align*}
		Based on the estimates above, using Lemmas \ref{lemma3.3} and \ref{lemma3.5}, we get \eqref{3-91}.
		We thus finish the proof of Lemma \ref{lemma3.11}.
	\end{proof}
	
	\begin{proof}[Proof of Proposition \ref{proposition3.2}]
		Choose $\delta_1 \leq \delta_0$ sufficiently small that all estimates in Lemmas \ref{lemma3.3}--\ref{lemma3.11} apply. If $N(t) \leq \delta_1$ for any $0 \leq t \leq T$, then estimate \eqref{3-12} follows immediately. Hence Proposition \ref{proposition3.2} is proved.
	\end{proof}
	
	\subsection{Global existence and convergence}
	\begin{proof}[Proof of Theorem \ref{theorem3.1}]
		Choose $\delta_{2} > 0$ sufficiently small that $\delta_2 = \min\{\delta_1/(2C_{*}), \delta_1/(C C_{*})\}$. By the standard continuation arguments associated with Propositions \ref{proposition3.1} and \ref{proposition3.2}, we obtain a unique global solution $(\phi, \psi) \in X(0, \infty)$ satisfying estimate \eqref{3-12} for all $t \in[0, \infty)$.
		
		Before proving \eqref{3-13}, we recall a classical result for later use; see \cite{MANKG}.
		\begin{lemma}\label{lemma3.10}
			If $F\in W^{1,1}(0,\infty)$, then $F(t)\to0$ as $t\to\infty$.
		\end{lemma}
		\begin{proof}
			Since $F'\in L^1$, $F$ has a finite limit at infinity.
			Since $F\in L^1$, the limit must be zero.
		\end{proof}
		
		To complete the proof of Theorem \ref{theorem3.1}, we exploit the global estimate \eqref{3-12}, which has already been established for all $t \in[0, \infty)$. The global estimate \eqref{3-12} implies
		\begin{equation}
			\|r\|_{w_{2}}^2 + \|r_{z}\|_{w_{2}}^2 + \|q\|_{w_{2}}^2 + \|q_{z}\|_{w_{2}}^2 + \int_0^{\infty} (\|r\|_{w_{1}}^2 + \|r_{z}\|_{w_{1}}^2 + \|q\|_{w_{1}}^2 + \|q_{z}\|_{1,w_{3}}^2) d \tau \leq C.
			\label{3-612}
		\end{equation}
		In particular, $\|r(\cdot, t)\|^2 \in L^1(0, \infty)$. It remains to verify that $\frac{\mathrm{d}}{\mathrm{d}t}\|r(\cdot, t)\|^2 \in L^1(0, \infty)$. By the Cauchy--Schwarz inequality and the identity $r_{t} = sr_{z} + q_{z}$ of \eqref{3-251},
		\begin{align*}
			\int_0^{\infty}\left|\frac{d}{d \tau} \|r(\cdot, t)\|^2\right| d \tau = \int_0^{\infty}\left|\frac{d}{d \tau} \int r^2 \right| d \tau = 2 \int_0^{\infty}\left|\int r r_{t} \right| d \tau \leq C \int_0^{\infty} \int (r^2 + r_{z}^2 + q_{z}^2) d \tau \leq C.
		\end{align*}
		Using Lemma \ref{lemma3.10}, we conclude that
		\begin{align*}
			\|r(\cdot, t)\|^2 \to 0, \quad \text { as } t \to +\infty.
		\end{align*}
		Hence, by \eqref{3-612}, for all $z \in \mathbb{R}$, it holds that
		\begin{align*}
			r^2(z, t) = & 2 \int_{-\infty}^z r r_{z}(y, t) d y
			\\
			\leq & 2\left(\int_{-\infty}^{\infty} r^2 d y\right)^{1 / 2}\left(\int_{-\infty}^{\infty} r_{z}^2 d y\right)^{1 / 2}
			\\
			\leq & C \left\|r(\cdot, t)\right\| \to 0, \, \text { as } t \to +\infty,
		\end{align*}
		namely,
		$$
		\sup _{z \in \mathbb{R}}|r(z, t)| \to 0, \quad \text { as } t \to + \infty.
		$$
		Recall from \eqref{4.7} that
		\begin{equation*}
			q_t = sq_z  - (p(V+r)-p(V))_{z} + \mu \left(\frac{q_z}{V+r} + \frac{sV_zr}{V(V+r)}\right)_z.
		\end{equation*}
		By \eqref{3-612} and Taylor's expansion, it is obvious that
		\begin{align*}
			\int_0^t \int (s q_z)^2 \leq C,
		\end{align*}
		and
		\begin{align*}
			\int_0^\infty \int [-(p(V+r)-p(V))_z]^2 = & \int_0^\infty \int [p'(V+r)r_z + V_z (p'(V+r)-p'(V))]^2
			\\
			\leq & C \int_0^\infty \int (r_{z}^{2}+ V_z^{2} r^{2})
			\\
			\leq & C \int_0^\infty \left(\int \frac{r_{z}^2}{V}+\int \frac{r^2}{V}\right)
			\\
			\leq & C (\|\phi_0\|_1^2 + \|\psi_0\|_1^2).
		\end{align*}
		Furthermore, by \eqref{3-16}, \eqref{4.711}, \eqref{3-612}, $\| r_{z}/V \|_{L^{\infty}} \leq C N(t)$ and Taylor's expansion, we derive
		\begin{align*}
			& \int_0^\infty \int \mu^{2} \left[\left(\frac{q_z}{V+r}\right)_z\right]^{2} + \int_0^{\infty} \int \mu^{2}  \left[\left(\frac{sV_zr}{V(V+r)}\right)_z\right]^2
			\\
			= & \mu^{2} \int_0^\infty \int \left(\frac{q_{zz}}{V + r} - \frac{q_{z}(V_{z} + r_{z})}{(V + r)^{2}}\right)^{2} + \mu^{2} \int_0^\infty \int \left[ -s \left(\frac{ V_{z z}}{V+r}-\frac{ V_{z z}}{V}\right) + s \left(\frac{V_z+r_{z}}{(V+r)^2}-\frac{V_z}{V^2}\right) V_z\right]^{2}
			\\
			\leq & C \int_0^\infty \int \left(\frac{q_{zz}^{2}}{(V+r)^{2}} + \frac{|q_{z} q_{zz} (V_{z} + r_{z})|}{(V+r)^{3}} + \frac{q_{z}^{2} (V_{z} + r_{z})^{2}}{(V+r)^{4}}\right)
			\\
			& + C \int_0^\infty \int \left(\frac{\left|V_{z z}\right|}{V^2}\left|r\right|+\frac{V_z^2}{V^3}\left|r\right|+\frac{\left|V_z\right|}{V^2}\left|r_{z}\right|\right)^2
			\\
			\leq & C \int_0^\infty \int \left(\frac{q_{zz}^{2}}{V^{2}} + \frac{|q_{z}q_{zz}|}{V} + \frac{|q_{z}q_{zz}r_{z}|}{V^{3}} + \frac{q_{z}^{2}}{V^{2}} + \frac{q_{z}^{2} r_{z}^{2}}{V^{4}} + V^{2} r^{2} + r_{z}^{2}\right)
			\\
			\leq & C \int_0^\infty \int \left(\frac{q_{zz}^{2}}{V^{2}} + \frac{q_{z}^{2}}{V^{2}} + r_{z}^{2} + r^{2}\right)
			\\
			\leq & C.
		\end{align*}
		Applying the above estimates, we obtain
		\begin{align*}
			\int_0^{\infty} \int q_{t}^{2} \leq C.
		\end{align*}
		Combining this estimate with \eqref{3-612}, we have
		\begin{align*}
			\int_0^{\infty} \left|\frac{\mathrm{d}}{\mathrm{d} \tau} \|q(\cdot, \tau)\|^2\right| \mathrm{d} \tau  = \int_0^{\infty} \left|\frac{\mathrm{d}}{\mathrm{d} \tau} \int q^2 \right| \mathrm{d} \tau = 2 \int_0^{\infty} \left|\int q q_{t}\right| \mathrm{d} \tau \leq C \int_0^{\infty} \int (q^2 + q_{t}^2) \mathrm{d} \tau \leq C.
		\end{align*}
		By Lemma \ref{lemma3.10}, it is obvious that
		\begin{align}
			\|q(\cdot, t)\|^2 \to 0, \, \text { as } t \to +\infty.
			\label{4-6}
		\end{align}
		Hence, for all $z \in \mathbb{R}$, from \eqref{3-612}, \eqref{4-6} and the Sobolev inequality,
		\begin{align*}
			q^2(z, t) = & 2 \int_{-\infty}^z q q_{z}(y, t) \mathrm{d} y
			\\
			\leq & C \left(\int_{-\infty}^{\infty} q^2 \mathrm{d} y\right)^{1 / 2} \left(\int_{-\infty}^{\infty} q_{z}^2 \mathrm{d} y\right)^{1 / 2}
			\\
			\leq & C \|q(\cdot, t)\| \to 0, \, \text { as } t \to +\infty,
		\end{align*}
		namely,
		$$
		\sup _{z \in \mathbb{R}}|q(z, t)| \to 0, \quad \text { as } t \to \infty.
		$$
		Now we complete the proof of Theorem \ref{theorem3.1}.
	\end{proof}

	\begin{proof}[Proof of Theorem \ref{theorem1.3}]
		First take $x_0=0$. By the definitions of $b_0$ and $B$, we have
		\[
		q_0=\psi_{0z}+\mu\left[\log(v_0/V)\right]_z=u_0-U.
		\]
		Thus the smallness assumption implies the hypotheses of
		Theorem \ref{theorem3.1}. We reconstruct
		\begin{equation*}
			v(z+st,t)=V(z)+r(z,t),\qquad
			u(z+st,t)=U(z)+q(z,t).
		\end{equation*}
		The identity $r_t=s r_z+q_z$, together with $U_z=-sV_z$, gives the mass equation. \eqref{4.7} gives
		\begin{align*}
			q_t = s q_z-(p(v)-p(V))_z + \mu \left[\left(\frac{u_{z}}{v}-\frac{U_{z}}{V}\right)\right]_z.
		\end{align*}
		Combining \eqref{1-111}, \eqref{1-5}, \eqref{1-10} and \eqref{2-15} with the preceding identity yields the momentum equation
		\begin{align*}
			u_t-s u_z+p(v)_z = \mu \left(\frac{u_z}{v}\right)_z.
		\end{align*}
		The regularity and uniqueness follow from the
		perturbation theorem and $V+r \geq V/2>0$.
		Convergence follows directly from \eqref{3-13}.
		Restoring the translation $x_0$ gives the statement of the theorem.
	\end{proof}

	\appendix
	\section{Local existence and the singular limit}\label{sec:local}
	We prove Proposition \ref{proposition3.1}. The key point is to make nonsingularization for the profile and the initial data.
	
	\subsection{Nonsingularized initial data}
	For $0<\varepsilon\leq1$, set
	\begin{equation*}
		V_{\varepsilon} = V+\varepsilon,\qquad
		w_i^\varepsilon=V_{\varepsilon}^{-i},\quad i=1,2,3.
	\end{equation*}
	Write $r_0=\phi_{0z}$ and $\eta_0=\log(1+r_0/V)$. Since $\|\eta_0\|\leq C\|r_0/V\|$, $\eta_{0z}=(q_0-\psi_{0z})/\mu$ and $\eta_{0zz}=(q_{0z}-\psi_{0zz})/\mu$, the assumptions of Proposition \ref{proposition3.1} imply
	\begin{equation*}
		\eta_0\in H^2,\qquad
		\|\eta_0\|_2\leq C N(0),
	\end{equation*}
	and $r_0/V \in H^1$, $\|r_0/V\|_{1} \leq C N(0)$. Moreover, differentiating $r_{0} = V(e^{\eta_0}-1)$ twice with respect to $z$, we deduce that
	\begin{equation}
		\begin{aligned}
			\|r_{0zz}/V\| = & \|V_{zz}(e^{\eta_0}-1)/V + 2 V_{z} e^{\eta_0} \eta_{0z}/V + e^{\eta_0} \eta_{0z}^{2} + e^{\eta_0} \eta_{0zz}\|
			\\
			\leq & C [\|(e^{\eta_0}-1)\| + \|e^{\eta_0}\|_{L^{\infty}} (\|\eta_{0z}\| + \|\eta_{0z}\|_{L^{\infty}} \|\eta_{0z}\| + \|\eta_{0zz}\|)]
			\\
			\leq & C N(0).
		\end{aligned}
		\label{5.511}
	\end{equation}
	Combining these estimates, we conclude that $r_0/V \in H^2$, $\|r_0/V\|_{2} \leq C N(0)$ and $\phi_0\in H^3$. Define
	\begin{equation}
		\eta_{0,\varepsilon}=\log(1+r_0/V_{\varepsilon}),\qquad
		\phi_{0,\varepsilon}=\phi_0,\qquad
		\psi_{0,\varepsilon}=\psi_0+\mu(\eta_0-\eta_{0,\varepsilon}).
		\label{5.6}
	\end{equation}
	It follows that
	\begin{equation}
		q_{0,\varepsilon} :=\psi_{0,\varepsilon z}
		+\mu\eta_{0,\varepsilon z}=q_0.
		\label{initial:qpreserved}
	\end{equation}
	For any fixed $\varepsilon>0$, it is clear that $\eta_{0,\varepsilon}, \phi_{0,\varepsilon}, \psi_{0,\varepsilon}, q_{0,\varepsilon} \to 0$ as $z \to \pm \infty$. By \eqref{3-151} and the assumptions of Proposition \ref{proposition3.1}, we know $|r_{0}/V| \leq 1/2$ and $V + r_{0} \geq V/2 > 0$, which implies
	\begin{equation}
		|r_{0}/V_{\varepsilon}| \leq 1/2, \quad V_{\varepsilon} + r_{0} \geq V_{\varepsilon}/2 > 0
		\label{5.61}
	\end{equation}
	Since $|\eta_{0,\varepsilon}| \leq C|r_{0}/V_{\varepsilon}|$ and $1/V_{\varepsilon} \leq 1/V$, we have
	\begin{equation}
		\|\eta_{0,\varepsilon}\| \leq C\|r_{0}/V_{\varepsilon}\| \leq C\|r_{0}/V\| \leq CN(0).
		\label{5.1}
	\end{equation}
	It is obvious that $V_{\varepsilon z} = V_{z}$ and $V_{\varepsilon zz} = V_{zz}$, the profile estimates established in \eqref{3-16} and Lemma \ref{lemma3.2} remain available for $V_{\varepsilon}$. Based on \eqref{4.711}, we obtain
	\begin{equation*}
		\|\eta_{0,\varepsilon z}\| = \left\|\frac{V_z+r_{0z}}{V_{\varepsilon}+r_{0}}-\frac{V_z}{V_{\varepsilon}}\right\| \leq C \left\|\frac{r_{0z}}{V_{\varepsilon}}\right\| + C \left\|\frac{V_z r_{0}}{V_{\varepsilon}^{2}}\right\| \leq C(\|r_{0z}/V\| + \|r_{0}\|) \leq CN(0),
	\end{equation*}
	and
	\begin{equation}
		\begin{aligned}
			& |\eta_{0,\varepsilon zz}|
			\\
			= & \left|\frac{V_{zz}+r_{0zz}}{V_{\varepsilon}+r_{0}} - \frac{(V_z+r_{0z})^2}{(V_{\varepsilon}+r_{0})^2} + \frac{V_z^2}{V_{\varepsilon}^2} - \frac{V_{zz}}{V_{\varepsilon}}\right|
			\\
			= & \left| \left(-\frac{1}{(V_{\varepsilon}+r_{0})^2} + \frac{1}{V_{\varepsilon}^2}\right)V_z^2 + \left(\frac{1}{V_{\varepsilon}+r_{0}} - \frac{1}{V_{\varepsilon}}\right)V_{zz} - \frac{r_{0z}^2}{(V_{\varepsilon}+r_{0})^2} - 2 \frac{V_z r_{0z}}{(V_{\varepsilon}+r_0)^2} + \frac{r_{0zz}}{V_{\varepsilon}+r_{0}} \right|
			\\
			\leq & C \left[\left(\frac{V_z^2}{V_{\varepsilon}^3} + \frac{|V_{zz}|}{V_{\varepsilon}^2}\right)|r_{0}| + \frac{r_{0z}^2}{V_{\varepsilon}^2} + \frac{|V_z||r_{0z}|}{V_{\varepsilon}^2} + \frac{|r_{0zz}|}{V_{\varepsilon}}\right].
		\end{aligned}
		\label{5.2}
	\end{equation}
	Thus, using \eqref{3-15} and \eqref{5.511}, it holds that
	\begin{equation}
		\begin{aligned}
			\|\eta_{0,\varepsilon zz}\| \leq C(\|r_{0}\| + \|r_{0z}/V\| + \|r_{0zz}/V\|) \leq CN(0).
		\end{aligned}
		\label{5.3}
	\end{equation}
	Combining \eqref{5.1}--\eqref{5.3}, we get
	\begin{equation}
		\|\eta_{0,\varepsilon}\|_2 \leq C N(0).
		\label{5.4}
	\end{equation}
	By $1/V_{\varepsilon} \leq 1/V$, \eqref{5.4} and the definitions of \eqref{5.6}, we have
	\begin{equation}
		\|\psi_{0,\varepsilon}\|_2 \leq C \|\psi_0\|_2 + C \|\eta_0\|_2 + C \|\eta_{0,\varepsilon}\|_{2} \leq C N(0),\quad
		\mathcal E_\varepsilon(0)\leq C\mathcal E(0),
		\label{initial:uniform}
	\end{equation}
	where
	\begin{equation*}
		\mathcal E_\varepsilon(t)=\|\phi^\varepsilon\|^2
		+\|\psi^\varepsilon\|_2^2
		+\|r^\varepsilon\|_{1,w_2^\varepsilon}^2
		+\|q^\varepsilon\|_{1,w_2^\varepsilon}^2.
	\end{equation*}
	For almost every $z$, since $V_{\varepsilon} = V+\varepsilon \to V$ as $\varepsilon\to0^+$, we obtain
	\begin{equation}
		\eta_{0,\varepsilon} \to \eta_0, \quad \eta_{0,\varepsilon z} \to \eta_{0z}, \quad \eta_{0,\varepsilon zz} \to \eta_{0zz}.
		\label{5.7}
	\end{equation}
	Using $|\eta_0| \leq C |r_0/V|$ and $\eta_{0z}=(q_0-\psi_{0z})/\mu$ yields
	\begin{equation}
		|\eta_0| + |\eta_{0z}| + |\eta_{0zz}| \leq C (|r_0/V| + |q_0| + |\psi_{0z}| + |q_{0z}| + |\psi_{0zz}|).
		\label{5.8}
	\end{equation}
	By virtue of \eqref{3-6112} and \eqref{5.2}, it is obvious that
	\begin{equation}
		|\eta_{0,\varepsilon}| + |\eta_{0,\varepsilon z}| + |\eta_{0,\varepsilon zz}| \leq C (|r_0/V| + |r_{0z}/V| + |r_{0zz}/V| + |q_0| + |\psi_{0z}| + |q_{0z}| + |\psi_{0zz}|).
		\label{5.9}
	\end{equation}
	Under the assumption of Proposition \ref{proposition3.1}, all functions on the right--hand sides of \eqref{5.8} and \eqref{5.9} belong to $L^{2}(\mathbb{R})$. Based on \eqref{5.7}, the dominated convergence theorem gives
	\begin{equation*}
		\|\eta_{0,\varepsilon}-\eta_0\|_{2} \to 0.
	\end{equation*}
	Consequently,
	\begin{equation}
		\eta_{0,\varepsilon}\to\eta_0\;\text{ in }H^2,\qquad
		\psi_{0,\varepsilon}\to\psi_0\;\text{ in }H^2.
		\label{initial:convergence}
	\end{equation}

	\subsection{Uniform existence time}
	We consider the regularized perturbation problem
	\begin{equation}
		\begin{cases}
			\phi_t^\varepsilon=s\phi_z^\varepsilon+\psi_z^\varepsilon
			+\mu\left(\dfrac{V_z+\phi_{zz}^\varepsilon}
			{V_{\varepsilon}+\phi_z^\varepsilon}-\dfrac{V_z}{V_{\varepsilon}}\right),\\
			\psi_t^\varepsilon=s\psi_z^\varepsilon
			-[p(V_{\varepsilon}+\phi_z^\varepsilon)-p(V_{\varepsilon})],\\
			(\phi^\varepsilon,\psi^\varepsilon)|_{t=0}
			=(\phi_{0,\varepsilon},\psi_{0,\varepsilon}).
		\end{cases}
		\label{A.1}
	\end{equation}
	Set
	\begin{equation}
		r^\varepsilon=\phi_z^\varepsilon,\qquad
		\eta^\varepsilon=\log(1+r^\varepsilon/V_{\varepsilon}),\qquad
		q^\varepsilon=\psi_z^\varepsilon+\mu\eta_z^\varepsilon
		=\phi_t^\varepsilon-s\phi_z^\varepsilon.
		\label{A.2}
	\end{equation}
	For each fixed $\varepsilon>0$, $V_{\varepsilon} \geq \varepsilon > 0$. By the standard fixed point theorems and $V_{\varepsilon} + r_{0} \geq V_{\varepsilon}/2 > 0$, we can prove local existence of solutions to the Cauchy problem \eqref{A.1}. This leads to the following local existence result.
	\begin{proposition}\label{propositionA-1}
		Suppose that $\phi_0 \in H^{3}(\mathbb{R})$ and $\psi_{0,\varepsilon} \in H^2(\mathbb{R})$. For any sufficiently small $\varepsilon_0 > 0$, there exists a positive constant $T_{\varepsilon}$ depending on $\varepsilon$ and $\varepsilon_0$ such that if $\|\phi_{0}\|_{3} + \|\psi_{0,\varepsilon}\|_{2} \leq \varepsilon_0$, then the problem \eqref{A.1} has a unique solution $(\phi^{\varepsilon}, \psi^{\varepsilon}) \in C([0,T_{\varepsilon}];H^3(\mathbb{R})) \times C([0,T_{\varepsilon}];H^2(\mathbb{R}))$ satisfying $\|\phi^{\varepsilon}\|_{3} + \|\psi^{\varepsilon}\|_{2} \leq C \varepsilon_{0}$ for any $0 \leq t \leq T_{\varepsilon}$.
	\end{proposition}
	Under assumptions of Proposition \ref{propositionA-1}, $V_{\varepsilon}+r^{\varepsilon} > 0$. One may first choose $J_n$ to be a standard Friedrichs mollifier and set $\phi_0^n=J_n*\phi_0$ and $\psi_{0,\varepsilon}^n=J_n*\psi_{0,\varepsilon}$. Then $\phi_0^n \to \phi_0$ in $H^3(\mathbb R)$ and $\psi_{0,\varepsilon}^n\to\psi_{0,\varepsilon}$ in $H^2(\mathbb R)$ as $n\to\infty$. By the Sobolev embedding and composition theorems, the induced initial values of $r^{\varepsilon}$ and $q^{\varepsilon}$ converge in $H^1_{w_2^\varepsilon}(\mathbb R)$. For sufficiently large $n$,  $V_{\varepsilon}+r_{0}^n \geq V_{\varepsilon}/2$ holds. Moreover, the previously established initial--energy estimate yields $\mathcal E_\varepsilon^n(0)\le C\mathcal E(0)$, where $C$ is independent of $\varepsilon$ and $n$. Since the subsequent energy estimates are uniform in both parameters, we pass to the limit $n\to\infty$ for each fixed $\varepsilon>0$. This yields a solution of the nonsingularized problem with the initial data $(\phi_{0},\psi_{0,\varepsilon})$.
	
	By \eqref{3-16}, we know that
	\begin{equation}
		|V_{\varepsilon z}| \leq C V^2, \quad |V_{\varepsilon zz}| \leq C V^3, \quad |V_{\varepsilon zzz}| \leq C V^4.
		\label{5.21}
	\end{equation}	
	Denote
	\begin{equation*}
		N_\varepsilon(t)
		:=\sup_{0\leq\tau\leq t}
		\mathcal E_\varepsilon(\tau)^{1/2}
		\leq\delta_*,
	\end{equation*}
	where \(\delta_*>0\) is a sufficiently small constant which is independent of $\varepsilon$ . Note that $\|r^{\varepsilon}/V_{\varepsilon}\|_{L^{\infty}} \leq C\|r^\varepsilon\|_{1,w_2^\varepsilon} \leq C N_\varepsilon(t)$. Therefore, $|r^\varepsilon|\leq V_{\varepsilon}/2$, $V_{\varepsilon}+r^\varepsilon \geq V_{\varepsilon}/2 \geq \varepsilon/2$. We define
	\begin{align*}
		d_\varepsilon=\frac{s p''(V_{\varepsilon})V_{\varepsilon z}}{p'(V_{\varepsilon})^2},\qquad
		\mathcal D_\varepsilon^{\rm loc}
		=\|r^\varepsilon\|_{2,w_1^\varepsilon}^2
		+\|q_z^\varepsilon\|_{1,w_3^\varepsilon}^2
		+\int d_\varepsilon
		[(\psi^\varepsilon)^2+(\psi_z^\varepsilon)^2
		+(\psi_{zz}^\varepsilon)^2].
	\end{align*}
	Repeating the proofs of Lemmas \ref{lemma3.3}, \ref{lemma3.5} and \ref{lemma3.11}, with $V_{\varepsilon}$ in place of $V$ gives
	\begin{equation}
		\|\phi^{\varepsilon}\|_{2}^2+\|\psi^{\varepsilon}\|_{2}^2+\int_0^t \|r^\varepsilon\|_{2,w_1^\varepsilon}^2 + \int_0^t \int d_\varepsilon \{(\psi^\varepsilon)^2+(\psi_z^\varepsilon)^2 +(\psi_{zz}^\varepsilon)^2\} \leq C (\|\phi_0\|_{2}^2 + \|\psi_{0\varepsilon}\|_{2}^2).
		\label{5.10}
	\end{equation}
	We next establish the weighted energy estimates of $q^{\varepsilon}$ and $q_{z}^{\varepsilon}$. Using the argument leading to \eqref{3-47} with $V_{\varepsilon}$ in place of $V$, we get
	\begin{align*}
		\|q^\varepsilon(t)\|_{w_2^\varepsilon}^2
		+\int_0^t\|q_z^\varepsilon\|_{w_3^\varepsilon}^2\,\mathrm{d}\tau \leq & C \|q_{0}\|_{w_{2}^{\varepsilon}}^{2} + C \int_{0}^{t} \int \frac{(q^\varepsilon)^{2}}{V_{\varepsilon}} + C \int_{0}^{t} \int \frac{(r^\varepsilon)^{2}}{V_{\varepsilon}}
		\\
		\leq & C\|q_0\|_{w_2^\varepsilon}^2
		+C\int_0^t\mathcal E_\varepsilon(\tau)\,\mathrm{d}\tau.
	\end{align*}
	By \eqref{4.13}, \eqref{5.21} and the above estimate, it follows that
	\begin{equation}
		\|q_z^\varepsilon(t)\|_{w_2^\varepsilon}^2
		+\int_0^t\|q_{zz}^\varepsilon\|_{w_3^\varepsilon}^2\,\mathrm{d}\tau \leq C\mathcal E_\varepsilon(0)
		+C\int_0^t\mathcal E_\varepsilon(\tau)\,\mathrm{d}\tau,
		\label{5.12}
	\end{equation}
	where $C > 0$ is independent of $t$ and $\varepsilon$. Next, we turn to the weighted $L^2$--estimates of $r^{\varepsilon}$ and $r_{z}^{\varepsilon}$. By Young's inequality, we have
	\begin{align*}
		\frac{1}{2} \frac{\mathrm{d}}{\mathrm{d}t} \int \frac{(r^{\varepsilon})^{2}}{V_{\varepsilon}^{2}} + \frac{1}{2} \frac{\mathrm{d}}{\mathrm{d}t} \int \frac{(r_{z}^{\varepsilon})^{2}}{V_{\varepsilon}^{2}} = & s \int \frac{r_{z}^{\varepsilon} r^{\varepsilon}}{V_{\varepsilon}^{2}} + s \int \frac{r_{zz}^{\varepsilon} r_{z}^{\varepsilon}}{V_{\varepsilon}^{2}} + \int \frac{q_{z}^{\varepsilon} r^{\varepsilon}}{V_{\varepsilon}^{2}} + \int  \frac{q_{zz}^{\varepsilon} r_{z}^{\varepsilon}}{V_{\varepsilon}^{2}}
		\\
		= & s \int \frac{ V_{z} (r^{\varepsilon})^{2}}{V_{\varepsilon}^{3}} + s \int \frac{ V_{z} (r_{z}^{\varepsilon})^{2}}{V_{\varepsilon}^{3}} + \int \frac{q_{z}^{\varepsilon} r^{\varepsilon}}{V_{\varepsilon}^{2}} + \int  \frac{q_{zz}^{\varepsilon} r_{z}^{\varepsilon}}{V_{\varepsilon}^{2}}
		\\
		\leq & C \int \frac{(q_{z}^\varepsilon)^{2}}{V_{\varepsilon}^{3}} + C \int \frac{(q_{zz}^\varepsilon)^{2}}{V_{\varepsilon}^{3}} + C \int \frac{(r^{\varepsilon})^{2}}{V_{\varepsilon}} + C \int \frac{(r_{z}^{\varepsilon})^{2}}{V_{\varepsilon}}.
		\label{A-6}
	\end{align*}
	Combining \eqref{5.10}--\eqref{5.12}, we then deduce that
	\begin{align*}
		\|r^\varepsilon\|_{w_2^\varepsilon}^2 + \|r_{z}^\varepsilon\|_{w_2^\varepsilon}^2 \leq C\mathcal E_\varepsilon(0)
		+C\int_0^t\mathcal E_\varepsilon(\tau)\,\mathrm{d}\tau.
	\end{align*}
	In summary, it holds that
	\begin{equation*}
		\mathcal E_\varepsilon(t)
		+\int_0^t\mathcal D_\varepsilon^{\rm loc}(\tau)\,\mathrm{d}\tau
		\leq C\mathcal E_\varepsilon(0)
		+C\int_0^t\mathcal E_\varepsilon(\tau)\,\mathrm{d}\tau.
	\end{equation*}
	Using the Gronwall inequality, we derive
	\begin{align}
		\mathcal E_\varepsilon(t) \leq C \mathcal E_\varepsilon(0) e^{C t},
		\label{A-9}
	\end{align}
	where $C > 0$ is independent of $t$ and $\varepsilon$. In particular, by virtue of \eqref{A.2}, we get
	\begin{align}
		\mu \eta_{zz}^{\varepsilon} = & \mu \left[\log \left(1 + \frac{r^{\varepsilon}}{V_{\varepsilon}}\right)\right]_{z z}
		\nonumber
		\\
		= &  \mu \left[\frac{V_{zz}+r_{zz}^{\varepsilon}}{V_{\varepsilon}+r^{\varepsilon}} - \frac{(V_{z}+r_{z}^{\varepsilon})^{2}}{(V_{\varepsilon}+r^{\varepsilon})^{2}} - \frac{V_{zz}}{V_{\varepsilon}} + \frac{V_{z}^{2}}{V_{\varepsilon}^{2}}\right]
		\nonumber
		\\
		= & \mu \left[\frac{r_{zz}^{\varepsilon}}{V_{\varepsilon}} + (V_{zz}+r_{zz}^{\varepsilon})\left(\frac{1}{V_{\varepsilon}+r^{\varepsilon}}  - \frac{1}{V_{\varepsilon}}\right) - \frac{(V_{z}+r_{z}^{\varepsilon})^{2}}{(V_{\varepsilon}+r^{\varepsilon})^{2}} + \frac{V_{z}^{2}}{V_{\varepsilon}^{2}}\right]
		\nonumber
		\\
		= & q_{z}^{\varepsilon} - \psi_{z z}^{\varepsilon}.
		\label{A-4111}
	\end{align}
	Under the assumptions of Proposition \ref{propositionA-1}, by \eqref{A-9} and \eqref{A-4111}, it holds that
	\begin{align*}
		r_{zz}^{\varepsilon}/V_{\varepsilon} \in L^2(\mathbb{R}), \quad \|r_{zz}^{\varepsilon}/V_{\varepsilon}\|^{2} \leq C \mathcal E_\varepsilon(0) e^{C t}.
	\end{align*}
	From \eqref{A-9} and $\mathcal E_\varepsilon(0)\leq C\mathcal E(0) \leq CN^{2}(0)$ in \eqref{initial:uniform}, we have
	\begin{align*}
		N_{\varepsilon}^{2}(t) \leq C \mathcal E_{\varepsilon}(0) e^{C t} \leq C N^{2}(0) e^{C t}.
	\end{align*}
	There exists $T_{0} > 0$ independent of $\varepsilon$ such that
	\begin{align*}
		N_{\varepsilon}^{2}(t) \leq C_{1} N^{2}(0) e^{C_{2} t}, \quad \hbox{for all } 0 < \varepsilon \leq 1,
	\end{align*}
	where $C_{1}, C_{2} > 0$ are constants independent of $t$ and $\varepsilon$ for any $t \in [0, T_{0}]$. If we choose $T_{0} = \frac{1}{C_{2}}$, we get
	\begin{align}
		N_{\varepsilon}(t) \leq C_{1}^{\frac{1}{2}} N(0) e^{\frac{1}{2}} =: C_{*} N(0),
		\label{A.14}
	\end{align}
	where $C_{*} = C_{1}^{\frac{1}{2}} e^{\frac{1}{2}}$ is a constant which is independent of $t$ and $\varepsilon$. Furthermore, as long as we choose $N(0) \leq \delta_0$ ($\delta_0$ is sufficiently small), $N_{\varepsilon}(t) \leq \delta_{*}/2$ holds on $[0,T_{0}]$.

	\subsection{Compactness and passage to the limit}
	For any $0 < \varepsilon \leq 1$, it follows from the above estimates that
	\begin{align}
		\begin{cases}
			\phi^{\varepsilon} \in L^{\infty}(0, T_{0}; H^{2}(\mathbb{R})), \quad \psi^{\varepsilon} \in L^{\infty}(0, T_{0}; H^{2}(\mathbb{R}));
			\\
			r^{\varepsilon} \in L^{\infty}(0, T_{0}; H_{w_{2}^{\varepsilon}}^{1}(\mathbb{R})), \quad q^{\varepsilon} \in L^{\infty}(0, T_{0}; H_{w_{2}^{\varepsilon}}^{1}(\mathbb{R})).
		\end{cases}
		\label{A.15}
	\end{align}
	Therefore, for $N(0) \leq \delta_{0}$, we have $1 + r/V_{\varepsilon} \geq 1 - CN(0) > 0$ where $C$ is independent of $\varepsilon$. Set $K_{M} := [-M, M] \subset \mathbb{R}$. $V > 0$ on $K_{M}$ ensures the existence of $c_K>0$ such that
	\begin{equation}
		V_{\varepsilon} \geq V \geq c_K > 0.
		\label{A.16}
	\end{equation}
	From \eqref{A.1}$_{2}$, \eqref{A.15}, $q^{\varepsilon} = \phi_{t}^{\varepsilon} - s \phi_{z}^{\varepsilon}$ and $\phi_t^\varepsilon=s r^\varepsilon+q^\varepsilon$, we have
	\begin{align*}
		\begin{cases}
			\| \phi_{t}^{\varepsilon} \| \leq \| q^{\varepsilon} \| + C \| r^{\varepsilon} \| \leq C, \quad \| \psi_t^{\varepsilon} \| \leq C \| \psi_z^{\varepsilon} \| + C \| r^{\varepsilon} \| \leq C,
			\\
			\| r_{t}^{\varepsilon} \| \leq \| q_{z}^{\varepsilon} \| + C \| r_{z}^{\varepsilon} \| \leq C, \quad \| \psi_{zt}^{\varepsilon} \| \leq C \| \psi_{zz}^{\varepsilon} \| + C \| r_{z}^{\varepsilon} \| + C |V_{z}| \| r^{\varepsilon} \| \leq C,
		\end{cases}
	\end{align*}
	for any $t \in [0, T_{0}]$. Namely, $\phi_{t}^{\varepsilon}, \psi_t^{\varepsilon}, r_{t}^{\varepsilon}, \psi_{zt}^{\varepsilon} \in  L^{\infty}(0, T_{0}; L^{2}(\mathbb{R}))$. Simon's compactness theorem, together with Sobolev embedding, interpolation and the Arzel\`a--Ascoli theorem, gives
	\begin{align*}
		\begin{cases}
			\phi^{\varepsilon} \to \phi \quad \text{strongly in } C([0, T_{0}]; H^{1}(K_{M})), \quad \text{as } \varepsilon \to 0,
			\\
			\psi^{\varepsilon} \to \psi \quad \text{strongly in } C([0, T_{0}]; H^{1}(K_{M})), \quad \text{as } \varepsilon \to 0.
		\end{cases}
	\end{align*}
	Since the compact interval $K_{M}$ is arbitrary, a standard diagonal argument allows us to extract a single subsequence satisfying
	\begin{align}
		\begin{cases}
			\phi^{\varepsilon} \to \phi \quad \text{strongly in } C([0, T_{0}]; H_{\mathrm{loc}}^{1}(\mathbb{R})), \quad \text{as } \varepsilon \to 0,
			\\
			\psi^{\varepsilon} \to \psi \quad \text{strongly in } C([0, T_{0}]; H_{\mathrm{loc}}^{1}(\mathbb{R})), \quad \text{as } \varepsilon \to 0.
		\end{cases}
		\label{A.17}
	\end{align}
	Moreover,
	\begin{align}
		\begin{cases}
			r_{z}^{\varepsilon} \rightharpoonup^* r_{z} \quad \text{in } L^{\infty}(0, T_{0}; L_{\mathrm{loc}}^{2}(\mathbb{R})), \quad \text{as } \varepsilon \to 0,
			\\
			\psi_{zz}^{\varepsilon} \rightharpoonup^* \psi_{zz} \quad \text{in } L^{\infty}(0, T_{0}; L_{\mathrm{loc}}^{2}(\mathbb{R})), \quad \text{as } \varepsilon \to 0,
		\end{cases}
		\label{A.18}
	\end{align}
	follows from \eqref{A.15}. Combining \eqref{initial:convergence}, \eqref{A.17}$_{2}$ and \eqref{A.18}$_{2}$, we have
	\begin{equation}
		\psi(\cdot,0) = \psi_0
		\quad\text{in }H_{\mathrm{loc}}^{1}(\mathbb{R}); \quad \psi_{zz}(\cdot,t)\rightharpoonup\psi_{0zz} \quad\text{weakly in }L_{\mathrm{loc}}^{2}(\mathbb{R}),
		\quad\text{as }t\to0^+.
		\label{A.21}
	\end{equation}
	The identity $V_{\varepsilon} = V + \varepsilon$ and \eqref{A.17} imply
	\begin{align*}
		V_{\varepsilon} + r^{\varepsilon} \to V + r \quad \text{strongly in } C([0, T_{0}]; L_{\mathrm{loc}}^{2}(\mathbb{R})), \quad \text{as } \varepsilon \to 0.
	\end{align*}
	Then, in view of \eqref{A.16}, we know $V_{\varepsilon} + r^{\varepsilon} \geq V_{\varepsilon}/2 \geq c_K/2 > 0$ and $V + r \geq V/2 \geq c_K/2 > 0$. Then we further obtain
	\begin{align}
		\frac{1}{V_{\varepsilon} + r^{\varepsilon}} \to \frac{1}{V + r} \quad \text{strongly in } C([0, T_{0}]; L_{\mathrm{loc}}^{2}(\mathbb{R})), \quad \text{as } \varepsilon \to 0.
		\label{A.19}
	\end{align}
	By virtue of \eqref{A.18}, \eqref{A.19} and $p \in C^4$, using the strong--weak convergence principle gives, for $\varepsilon \to 0$,
	\begin{align}
		\begin{cases}
			\frac{r_{z}^{\varepsilon}}{V_{\varepsilon} + r^{\varepsilon}} \rightharpoonup^* \frac{r_{z}}{V + r} \quad \text{in } L^{\infty}(0, T_{0}; L_{\mathrm{loc}}^{2}(\mathbb{R})),
			\\
			\frac{V_{z}}{V_{\varepsilon} + r^{\varepsilon}} \to \frac{V_{z}}{V +r},\; \frac{V_{z}}{V_{\varepsilon}} \to \frac{V_{z}}{V},\; p(V_{\varepsilon} + r^{\varepsilon}) \to p(V + r),\; p(V_{\varepsilon}) \to p(V) \; \text{in } C([0, T_{0}]; L_{\mathrm{loc}}^{2}(\mathbb{R})).
		\end{cases}
		\label{A.20}
	\end{align}
	Let $g \in C_{c}^{\infty}(\mathbb{R} \times (0, T_{0}))$ be arbitrary. We can make $\operatorname{supp}g \subset K_M \times (0,T_0)$. Hence, all the strong and weak convergence results established above on \(K_M\) are applicable. Multiplying \eqref{A.1}$_{1}$ by $g$, integrating it over $\mathbb{R} \times (0,T_0)$ and integrating by parts with respect to $t$, we obtain
	\begin{align*}
		\int_0^{T_0} \int \phi_{t}^\varepsilon g = -\int_0^{T_0} \int \phi^\varepsilon g_t = \int_0^{T_0} \int \left[ sr^\varepsilon+\psi_z^\varepsilon + \mu \left( \frac{V_z+r_{z}^\varepsilon}{V_{\varepsilon}+r^\varepsilon} - \frac{V_z}{V_{\varepsilon}} \right) \right] g.
	\end{align*}
	Letting $\varepsilon \to 0$ and using \eqref{A.17} and \eqref{A.20}, we infer that
	\begin{align*}
		-\int_0^{T_0} \int \phi g_t = \int_0^{T_0} \int \left[ sr + \psi_z + \mu \left( \frac{V_z+r_{z}}{V + r} - \frac{V_z}{V} \right) \right] g.
	\end{align*}
	It follows that
	\begin{align}
		\phi_t = s\phi_z + \psi_z + \mu \left(\frac{V_z+\phi_{zz}}{V+\phi_z} - \frac{V_z}{V}\right)
		\label{A.22}
	\end{align}
	holds in the sense of distributions in $\mathbb{R} \times (0,T_0)$. The second equation follows in the same manner. Thus,
	\begin{align}
		\psi_t = s\psi_z - [p(V+\phi_z)-p(V)]
		\label{A.231}
	\end{align}
	holds in the sense of distributions in $\mathbb{R} \times (0,T_0)$. The initial values of $\phi$ and $\psi$ follow from $\phi_{0,\varepsilon} = \phi_{0}$, $\psi_{0,\varepsilon} \to \psi_{0}$ in $H_{\mathrm{loc}}^{2}(\mathbb{R})$, the local convergence \eqref{A.17}--\eqref{A.18} and \eqref{A.21}. Consequently, by virtue of \eqref{A.22} and \eqref{A.231}, $(\phi,\psi)$ satisfies \eqref{3-611} in $\mathbb{R} \times (0,T_0)$ in the sense of distributions.
	
	By the weak lower semicontinuity of the $L^2$--norm and \eqref{A.14}, we obtain
	\begin{equation}
		\sup_{t\leq T_0}\mathcal E(t)
		+\int_0^{T_0}\left\{
		\|r\|_{2,w_1}^2+\|q_z\|_{1,w_3}^2
		+\int d(V)(\psi^2+\psi_z^2+\psi_{zz}^2)\right\} \mathrm{d}t
		\leq C\mathcal E(0).
		\label{A.24}
	\end{equation}
	
	\subsection{Strong time continuity}
	Based on \eqref{A.24}, we get
	\begin{equation}
		\begin{cases}
			\phi \in L^{\infty}(0, T_{0}; L^{2}(\mathbb{R})), \quad \psi \in L^{\infty}(0, T_{0}; H^{2}(\mathbb{R}));
			\\
			r \in L^{\infty}(0, T_{0}; H_{w_{2}}^{1}(\mathbb{R})), \quad q \in L^{\infty}(0, T_{0}; H_{w_{2}}^{1}(\mathbb{R})).
		\end{cases}
		\label{A.25}
	\end{equation}
	By virtue of $\phi_{t} = q + s r$, \eqref{3-611}$_{2}$ and \eqref{A.25}, it holds that
	\begin{align*}
		\begin{cases}
			\| \phi_{t} \|^{2} \leq C \| q \|^{2} + C \| r \|^{2} \leq C, \quad \| \psi_t \|^{2} \leq C \| \psi_z \|^{2} + C \| r \|^{2} \leq C,
			\\
			\| \frac{r_{t}}{V} \|^{2} \leq C \| \frac{q_{z}}{V} \|^{2} + C \| \frac{r_{z}}{V} \|^{2} \leq C, \quad \| \psi_{zt} \|^{2} \leq C \| \psi_{zz} \|^{2} + C \| r_{z} \|^{2} + C |V_{z}|^{2} \| r \|^{2} \leq C.
		\end{cases}
	\end{align*}
	Consequently,
	\begin{align}
		\begin{cases}
			\phi_{t} \in L^{\infty}(0, T_{0}; L^{2}(\mathbb{R})), \quad \psi_{t} \in L^{\infty}(0, T_{0}; L^{2}(\mathbb{R}));
			\\
			r_{t} \in L^{\infty}(0, T_{0}; L_{w_{2}}^{2}(\mathbb{R})), \quad \psi_{zt} \in L^{\infty}(0, T_{0}; L^{2}(\mathbb{R})).
		\end{cases}
		\label{A.26}
	\end{align}
	Combining \eqref{A.25} and \eqref{A.26}, we have
	\begin{align}
		\phi, \, \psi, \, \psi_{z} \in C([0, T_{0}]; L^{2}(\mathbb{R})), \quad r \in C([0, T_{0}]; L_{w_{2}}^{2}(\mathbb{R})).
		\label{A.27}
	\end{align}
	To obtain the desired continuity, we denote
	\begin{equation}
		\begin{aligned}
			\mathcal{R}_{1}(t) := \frac{\mathrm{d}}{\mathrm{d}t} \| q \|_{w_{2}}^2, \quad \mathcal{R}_{2}(t) := \frac{\mathrm{d}}{\mathrm{d}t} \| q_{z} \|_{w_{2}}^2, \quad \mathcal{R}_{3}(t) := \frac{\mathrm{d}}{\mathrm{d}t} \| \psi_{zz} \|^2, \quad \mathcal{R}_{4}(t) := \frac{\mathrm{d}}{\mathrm{d}t} \| r_{z} \|_{w_{2}}^2.
		\end{aligned}
		\label{A.28}
	\end{equation}
	Note that $\int$\eqref{3-4611}, $\int$\eqref{4.9}, $\int$\eqref{3-411} and
	\begin{equation*}
		\frac{1}{2} \int \frac{\mathrm{d}}{\mathrm{d}t} \frac{r_z^2}{V^2} = \int \frac{q_{zz}r_z}{V^2} + s \int \frac{V_zr_z^2}{V^3} + \frac{s}{2} \int \left(\frac{r_{z}^{2}}{V^{2}}\right)_{z}
	\end{equation*}
	provide specific expressions for \eqref{A.28}, respectively. Using \eqref{A.24} and a method similar to that in Lemmas \ref{lemma3.6}--\ref{lemma3.11}, we clearly derive that
	\begin{equation*}
		\int_{0}^{T_{0}} (|\mathcal{R}_{1}(t)| + |\mathcal{R}_{2}(t)| + |\mathcal{R}_{3}(t)| + |\mathcal{R}_{4}(t)| ) \leq C\mathcal E(0),
	\end{equation*}
	namely, $\mathcal{R}_{i} \in L^1(0,T_0)$, $i=1,2,3,4$. It is obvious that $\| r_{z} \|_{w_{2}}^2, \, \| \psi_{zz} \|^2, \, \| q \|_{w_{2}}^2, \, \| q_{z} \|_{w_{2}}^2 \in L^1(0,T_0)$ holds by \eqref{A.24}. Therefore, it is clear to know that
	\begin{align*}
		\| r_{z} \|_{w_{2}}^2, \, \| \psi_{zz} \|^2, \, \| q \|_{w_{2}}^2, \, \| q_{z} \|_{w_{2}}^2 \in W^{1,1}(0,T_0).
	\end{align*}
	In the one--dimensional case, the Sobolev space $W^{1,1}(0,T_0)$ embeds continuously into $C([0,T_0])$. We conclude that
	\begin{equation}
		\| r_{z} \|_{w_{2}}^2, \, \| \psi_{zz} \|^2, \, \| q \|_{w_{2}}^2, \, \| q_{z} \|_{w_{2}}^2 \in W^{1,1}(0,T_0) \hookrightarrow C([0,T_0]).
		\label{A.291}
	\end{equation}
	
	Next, we establish weak continuity of $r_{z}/V$, $\psi_{zz}$, $q/V$ and $q_{z}/V$. We present the argument for $q/V$, as the remaining cases can be treated in a similar manner.
	
	We hope to prove that $q/V \in L^2(\mathbb{R})$ for every $t \in [0,T_0]$. By \eqref{A.24}, \eqref{A.25} and \eqref{A.27}, we know that there exists a set of full measure $E \subset (0,T_0)$ such that, for every $t \in E$, $q=(\psi+\mu\eta)_{z}$ satisfies
	\begin{equation}
		\psi + \mu \eta \in C([0,T_0];L^2(\mathbb R)), \quad \frac{q}{V}\in L^\infty(0,T_0;L^2(\mathbb R)), \quad \left\|\frac{q}{V}(t)\right\|^{2} \leq C \mathcal{E}(0), \quad \hbox{for } t \in E.
		\label{5.22}
	\end{equation}
	Let $t_0 \in [0,T_0]$ be arbitrary. We can choose a sequence $t_n\in E$ with $t_n\to t_0$ as $n \to +\infty$. For each $\varphi \in C_c^\infty(\mathbb R)$, integration by parts gives
	\begin{equation*}
		\int \frac{q}{V}(z,t_n) \varphi(z) = - \int (\psi+\mu\eta)(z,t_n) \left(\frac{\varphi(z)}{V}\right)_{z}.
	\end{equation*}
	Using \eqref{5.22}, we know that
	\begin{equation}
		\|(\psi+\mu\eta)(t_n)-(\psi+\mu\eta)(t_0)\| \to 0, \quad \hbox{as } t_n\to t_0.
		\label{5.25}
	\end{equation}
	Since there exists a constant $c_{1}>0$ such that $V \geq c_{1} > 0$ for all $z \in \operatorname{supp}\varphi$, we have
	\begin{equation}
		\left\| \left(\frac{\varphi}{V}\right)_{z} \right\| < +\infty.
		\label{5.24}
	\end{equation}
	By \eqref{5.25}, \eqref{5.24} and the Cauchy--Schwarz inequality, as $t_n\to t_0$, it follows that
	\begin{align*}
		\left|\int [(\psi+\mu\eta)(z,t_n) - (\psi+\mu\eta)(z,t_0)] \left(\frac{\varphi(z)}{V}\right)_{z}\right| \leq \|(\psi+\mu\eta)(t_n)-(\psi+\mu\eta)(t_0)\| \left\| \left(\frac{\varphi}{V}\right)_{z} \right\| \to 0.
	\end{align*}
	Consequently,
	\begin{align*}
		\lim_{n \to +\infty} \int \frac{q}{V}(z,t_n) \varphi(z) = - \int (\psi+\mu\eta)(z,t_0) \left(\frac{\varphi(z)}{V}\right)_{z}.
	\end{align*}
	Combining \eqref{5.22} and $t_{n} \in E$, we obtain
	\begin{equation*}
		\left|- \int (\psi+\mu\eta)(z,t_0) \left(\frac{\varphi(z)}{V}\right)_{z}\right| \leq \limsup_{n\to\infty} \left\|\frac{q}{V}(t_{n})\right\| \|\varphi\| \leq C \mathcal{E}^{\frac{1}{2}} (0) \|\varphi\|.
	\end{equation*}
	Since $C_c^\infty(\mathbb R)$ is dense in $L^2(\mathbb R)$, by the Riesz representation theorem, we know that there exists a unique function $g(z, t_0)\in L^2(\mathbb R)$ such that
	\begin{equation}
		\int g(z,t_0) \varphi(z) = - \int (\psi+\mu\eta)(z,t_0) \left(\frac{\varphi(z)}{V}\right)_{z},
		\label{A.40}
	\end{equation}
	for every $\varphi\in C_c^\infty(\mathbb R)$, and $\|g(t_{0})\| \leq C \mathcal{E}^{\frac{1}{2}} (0)$.
	If $t_0\in E$, $q(z,t_0)/V$ satisfies the same identity as in
	\eqref{A.40}. Hence, by the uniqueness in the Riesz representation
	theorem, we define $\frac{q}{V}(z, t_0) := g(z,t_0)$. Since $t_0\in[0,T_0]$ is arbitrary, for every $t \in [0,T_{0}]$, it holds that
	\begin{equation*}
		\frac{q}{V}(z,t) \in L^{2}(\mathbb{R}), \quad \sup_{t\in[0,T_0]} \left\|\frac{q}{V}(t)\right\| \leq C \mathcal{E}^{\frac{1}{2}} (0).
	\end{equation*}
	Thus, we have extended the definition of $q/V$ from almost every $t\in(0,T_0)$ to every $t\in[0,T_0]$.
	
	By the preceding construction, for every $t_{0} \in [0, T_{0}]$, it follows that
	\begin{equation}
		\frac{q}{V}(z,t_{0}) \in L^{2}(\mathbb{R}), \quad \left\|\frac{q}{V}(t_{0})\right\|^{2} \leq C \mathcal{E}(0).
		\label{A.29}
	\end{equation}
	Since $C_{c}^{\infty}(\mathbb{R})$ is dense in $L^{2}(\mathbb{R})$, for any $\delta>0$, there exist $\varphi_{\delta} \in C_{c}^{\infty}(\mathbb{R})$ and a constant $c_{\delta}>0$ such that $V \geq c_{\delta} > 0$ for all $z \in \operatorname{supp}\varphi_\delta$, and
	\begin{equation}
		\left\| \frac{\varphi_{\delta z}}{V} \right\| < +\infty, \quad \left\| \frac{q}{V}(t_{0}) - \varphi_{\delta}(z)\right\| < \delta.
		\label{A.30}
	\end{equation}
	For any $t \in [0,T_{0}]$, by $q = \psi_{z} + \mu \eta_{z}$, we have
	\begin{align*}
		& \left| \int \left(\frac{q}{V}(z,t) - \frac{q}{V}(z,t_{0})\right) \frac{q}{V}(z,t_{0}) \right|
		\\
		\leq & \left| \int \left(\frac{q}{V}(z,t) - \frac{q}{V}(z,t_{0})\right) \varphi_{\delta}(z) \right| + \left| \int \left(\frac{q}{V}(z,t) - \frac{q}{V}(z,t_{0})\right) \left(\frac{q}{V}(z,t_{0}) - \varphi_{\delta}(z)\right) \right|
		\\
		\leq & \left| \int [(\psi + \mu \eta)(z,t) - (\psi + \mu \eta)(z,t_{0})] \frac{\varphi_{\delta z}(z)}{V} \right| + \left| \int  [(\psi + \mu \eta)(z,t) - (\psi + \mu \eta)(z,t_{0})] \frac{V_{z} \varphi_{\delta}(z)}{V^{2}} \right|
		\\
		& + \left| \int \left(\frac{q}{V}(z,t) - \frac{q}{V}(z,t_{0})\right) \left(\frac{q}{V}(z,t_{0}) - \varphi_{\delta}(z)\right) \right|
		\\
		\leq & C\| (\psi + \mu \eta)(t) - (\psi + \mu \eta)(t_{0})\| \left\| \frac{\varphi_{\delta z}}{V} \right\| + C\| (\psi + \mu \eta)(t) - (\psi + \mu \eta)(t_{0})\| \| \varphi_{\delta} \|
		\\
		& + \left\| \frac{q}{V}(t) - \frac{q}{V}(t_{0}) \right\| \left\| \frac{q}{V}(t_{0}) - \varphi_{\delta} \right\|.
	\end{align*}
	For each fixed $\delta$, as $\eta = \log(1+r/V)$, $\psi \in C([0, T_{0}]; L^{2}(\mathbb{R}))$ and $r \in C([0, T_{0}]; L_{w_{2}}^{2}(\mathbb{R}))$, we obtain
	\begin{equation}
		\lim_{t\to t_0} \| (\psi + \mu \eta)(t) - (\psi + \mu \eta)(t_{0})\| \to 0.
		\label{A.32}
	\end{equation}
	Keeping $\delta > 0$ fixed, we get from \eqref{A.30} and \eqref{A.32} that
	\begin{equation*}
		\begin{aligned}
			\limsup_{t\to t_0} \left| \int \left(\frac{q}{V}(z,t) - \frac{q}{V}(z,t_{0})\right) \frac{q}{V}(z,t_{0}) \right| \leq \left(\left\| \frac{q}{V}(t) \right\| + \left\| \frac{q}{V}(t_{0}) \right\|\right) \left\| \frac{q}{V}(t_{0}) - \varphi_{\delta} \right\|
			\leq C \mathcal{E}^{\frac{1}{2}} (0) \delta.
		\end{aligned}
	\end{equation*}
	Then letting $\delta \to 0$, we derive that
	\begin{equation}
		\begin{aligned}
			\lim_{t\to t_0} \int \frac{q}{V}(z,t) \frac{q}{V}(z,t_{0}) = \int \left(\frac{q}{V}(z,t_{0})\right)^{2}.
		\end{aligned}
		\label{A.33}
	\end{equation}
	Combining \eqref{A.33} with the continuity of its $L^{2}$--norm in \eqref{A.291}, we get
	\begin{equation}
		\lim_{t\to t_0} \left\| \frac{q}{V}(t) - \frac{q}{V}(t_{0}) \right\|^{2} = \lim_{t\to t_0} \left(\left\| \frac{q}{V}(t) \right\|^{2} + \left\| \frac{q}{V}(t_{0}) \right\|^{2} - 2  \int \frac{q}{V}(z,t) \frac{q}{V}(z,t_{0}) \right) = 0.
		\label{A.35}
	\end{equation}
	Consequently, for every fixed $t_{0} \in [0,T_{0}]$, we arrive at
	\begin{equation}
		\frac{q}{V}(z,t) \to \frac{q}{V}(z,t_{0}) \quad \text{strongly in } L^{2}(\mathbb{R}), \quad \text{as } t\to t_0.
		\label{A.36}
	\end{equation}
	Finally, we have
	\begin{equation}
		q \in C([0, T_{0}]; L_{w_{2}}^{2}(\mathbb{R})).
		\label{A.34}
	\end{equation}
	Similarly, we conclude that
	\begin{align}
		\psi_{zz} \in C([0, T_{0}]; L^{2}(\mathbb{R})); \quad r_{z}, \, q_{z} \in C([0, T_{0}]; L_{w_{2}}^{2}(\mathbb{R})).
		\label{A.31}
	\end{align}
	By virtue of \eqref{A.27}, \eqref{A.34} and \eqref{A.31}, we obtain
	\begin{equation}
		(\phi, \psi) \in X(0, T_{0}).
		\label{5.32}
	\end{equation}
	For each fixed \(\varepsilon>0\), recalling $J_{n}$ as the Friedrichs mollifier, we can verify that the method previously used to prove the strong continuity of \eqref{5.32} in the time direction still holds.
	
	\subsection{Uniqueness}
	Let $(\phi_i,\psi_i)$, $i=1,2$, be two solutions of \eqref{3-611} and
	\begin{align*}
		M=\phi_1-\phi_2,\qquad P=\psi_1-\psi_2,\qquad
		v_i=V+\phi_{iz}.
	\end{align*}
	$(M,P)$ satisfies
	\begin{equation}
		\begin{cases}
			M_t=sM_z+P_z+\mu[\log(v_1/v_2)]_z,\\
			P_t=sP_z-[p(v_1)-p(v_2)].
		\end{cases}
		\label{5.33}
	\end{equation}
	Since $v_i$ are comparable to $V$,
	\begin{align*}
		M_z\log(v_1/v_2)\geq cM_z^2/V,\qquad
		|p(v_1)-p(v_2)|\leq C|M_z|.
	\end{align*}
	Multiplying \eqref{5.33}$_{1}$ by $M$ and \eqref{5.33}$_{2}$ by $P$, respectively, adding the resulting equations and integrating over $\mathbb{R} \times [0,t]$ gives
	\begin{align*}
		\frac{\mathrm{d}}{\mathrm{d}t}(\|M\|^2+\|P\|^2)
		+c\int\frac{M_z^2}{V}
		\leq C(\|M\|^2+\|P\|^2).
	\end{align*}
	Since $M(\cdot,0) = P(\cdot,0) = 0$, Gronwall's inequality implies $M \equiv P \equiv 0$. Hence the local solution is unique. This completes the proof of Proposition \ref{proposition3.1}.

	\subsection*{Acknowledgments}
	The research of J. Li was partially supported by the National Natural Science Foundation of China (No. 12371216). The research of C. Liu was partially supported by the Science and Technology Development Plan Project of Jilin Province, China (No. 20250102009JC). The research of M. Mei was partially supported by NSERC of Canada grant RGPIN 2022-03374 and the National Natural Science Foundation of China (No. W2431005). 
	
	\
	
	\noindent {\bf Declaration}. This is an original study, no AI is applied for mathematical ideas and proofs.

\end{document}